\documentclass[11pt]{article}
\usepackage{latexsym,amssymb,amsmath,amsfonts,amsthm,mathrsfs}
\usepackage{graphicx}
\usepackage{dcolumn}
\usepackage{bm}
\usepackage{color}
\usepackage{epsfig}
\usepackage{hyperref}
\usepackage{multirow}
\usepackage{tikz}
\usepackage{booktabs}
\usetikzlibrary{shapes.geometric}

\makeatletter
\newif\ifmsbmloaded@
\def\loadmsbm{\msbmloaded@true
  \font\tenmsb=msbm10 scaled 1\@ptsize00
  \font\sevenmsb=msbm7 scaled 1\@ptsize00
  \font\fivemsb=msbm5 scaled 1\@ptsize00
  \alloc@8\fam\chardef\sixt@@n\msbfam
  \textfont\msbfam=\tenmsb
  \scriptfont\msbfam=\sevenmsb
  \scriptscriptfont\msbfam=\fivemsb
  }
\def\nonmatherr@#1{\errmessage%
{LateX error: \string#1\space allowed only in math mode}}
\def\Bbb{\relax\ifmmode\expandafter\Bbb@\else
  \expandafter\nonmatherr@\expandafter\Bbb\fi}
\def\Bbb@#1{{\Bbb@@{#1}}}
\def\Bbb@@#1{\fam\msbfam\relax#1}
\@addtoreset{equation}{section}

\makeatother \loadmsbm

\def\u1{u_1}
\def\b1{b_1}

\newcommand{\beq}{\begin{equation}}
\newcommand{\eeq}{\end{equation}}
\newcommand{\ben}{\begin{eqnarray}}
\newcommand{\een}{\end{eqnarray}}
\newcommand{\beno}{\begin{eqnarray*}}
\newcommand{\eeno}{\end{eqnarray*}}

\newtheorem{thm}{Theorem}[section]
\newtheorem{lem}{Lemma}[section]

\newtheorem{prop}{Proposition}[section]
\newtheorem{cor}{Corollary}[section]
\newtheorem{defin}{Definition}[section]
\newtheorem{que}{Question}[section]
\begin{document}
\title{\bf  Global Well-Posedness for NLS with a Class of $H^s$-Supercritical Data  }

\footnotesize
\author{\bf Jinsheng Han \\
{\it \small Department of Mathematics, Sichuan Normal University, Chengdu 610068, P. R. China.}\\
   \bf  Baoxiang Wang  \\
    {\it \small LMAM, School of Mathematics, Peking University, Beijing 100871, P. R. China.}\\
    {\small (Emails: pkuhan@math.scnu.edu.cn, wbx@pku.edu.cn)}     }
\date{}
 \maketitle

\renewcommand{\thefootnote}{\fnsymbol{footnote}}

\vspace{-1.2in} \vspace{.9in} \vspace{0.2cm} {\bf Abstract.}  We study the Cauchy problem for NLS with a class of $H^s$-super-critical data
\begin{align}
& {\rm i}u_t +\Delta u+ \lambda |u|^{2\kappa} u =0,  \quad u(0)=u_0  \label{NLSabstract}
\end{align}
and  show that \eqref{NLSabstract} is globally well-posed and scattering in $\alpha$-modulation spaces $M^{s,\alpha}_{2,1}$ ($\alpha\in [0,1), \ s> d\alpha/2-\alpha/\kappa$, $\kappa\in \mathbb{N}$ and $\kappa \geq 2/d$) for the sufficiently small data. Moreover, NLS is ill-posed in $M^{s,\alpha}_{2,1}$ if $s< d\alpha/2-\alpha/\kappa$. In particular,  we obtain a class of initial data $u_0$ satisfying for any $M\gg 1$,
\begin{align}
\|u_0\|_2 \sim   M^{1/\kappa-d/2 },   \ \      \|u_0\|_\infty \ =\infty ,   \ \ \|u_0\|_{M^{s,\alpha}_{2,1}}     \geq  M^{(1-\alpha)/\kappa}, \ \ \  \|u_0\|_{B^{s(\kappa)}_{2,\infty}} =\infty     \nonumber
\end{align}
such that NLS is globally well-posed in $M^{s,\alpha}_{2,1}$ if $\kappa>2/d, \ \alpha\in [0,1)\  d\alpha/2-\alpha/\kappa <s < s(\kappa):= d/2-1/\kappa$. Such a kind of data are super-critical in $H^{s(\kappa)}$ and have infinite amplitude.

\vspace{0.2cm} {\bf Key words.} $\alpha$-modulation space, NLS, Global well-posedness, super-critical data in $H^s$.

 \vspace{0.2cm}


\normalsize

\section{Introduction}

In this paper we study the Cauchy problem for the nonlinear
Schr\"odinger equation (NLS):
\begin{align}
& {\rm i}u_t +\Delta u+ \lambda |u|^{2\kappa} u =0,  \quad u(0)=u_0, \label{NLS}
\end{align}
where $u(t,x)$ is a complex valued function of $(t,x) \in \mathbb{R}
\times \mathbb{R}^d$, ${\rm i}=\sqrt{-1}$, $u_t=\partial/\partial t,
$   $\Delta= \partial^2/\partial^2
x_1 +... + \partial^2/\partial^2 x_d$, $ \kappa \in \mathbb{N}$, $\lambda\in \mathbb{R}$, $u_0$ is a complex valued
functions of $x\in \mathbb{R}^n$. $\lambda<0$ ($\lambda>0$) is said to be the defocusing (focusing) case for NLS. The solutions of NLS satisfy the conservations of mass and energy:
\begin{align}
M(u(t)) & =  \|u(t)\|^2_{2} = M(u_0), \label{mass} \\
E(u(t)) & =  \|\nabla u (t)\|_2 - \frac{2\lambda}{1+\kappa} \|u(t)\|^{2+2\kappa}_{2+2\kappa}= E(u_0).  \label{energy}
\end{align}
It is well known that \eqref{NLS} is invariant under the scaling
\begin{align}
  u_{\sigma} (t,x) := \sigma^{1/\kappa} u(\sigma^2t, \sigma x), \ \ \sigma>0,  \label{scaling}
\end{align}
 which means that if $u$ solves \eqref{NLS}, so does $u_{\sigma}$ with initial data $\sigma^{1/\kappa} u_0(\sigma \cdot)$. Denote $s(\kappa) =d/2-1/\kappa$. $\dot H^{s(\kappa)}$  is  a scaling invariant Sobolev space for NLS, namely, $\|u_{\sigma}\|_{\dot H^s} = \|u \|_{\dot H^s}$ for all $\sigma>0$ if $s=s(\kappa)$, which is said to be a critical Sobolev space for NLS. If $s>s(\kappa)$ ($s<s(\kappa)$), then $H^s$ is said to be subcritical (supercritical) Sobolev spaces.  It is known that NLS is local well posed in all subcritical spaces $H^s$ with $s> s(\kappa)$, and globally well posed in critical space $H^{s(\kappa)}$ if initial data are sufficiently small; cf. Cazenave and Weissler \cite{CaWe90} and Nakamura and Ozawa \cite{NaOz98} for exponential growth nonlinearity. Moreover, NLS is ill posed in all super-critical spaces $H^s$ with $s<s(\kappa)$, cf. \cite{BiKePoSvVe96}.

In the $H^1$-subcritical and $H^1$-critical cases, a large amount of work has been devoted to the study of the global well-posedness and the existence of scattering operators of NLS (cf. \cite{Bo99,CoKeStTaTa08,DuThKeMe15,Na99,RyVi07,St81,Ts85} and references therein), where the conservations of mass and energy play important roles. However, for the  $H^1$-supercritical cases, the nonlinearity is out of the control of the energy space and up to now, it is not very clear for the mechanism of the  well/ill posedness of NLS. For the defocusing energy-supercritical NLS ($\lambda<0, \ \kappa>2/(d-2)$), the global well-posedness for large data has been open for many years.  In the general $H^s$ super-critical cases, we have a more delicate question:

 \begin{que} \label{question}
  Are there any initial data out of the critical Sobolev spaces $ \dot H^{s(\kappa)}$ (or more general critical Besov spaces $ \dot B^{s(\kappa)}_{2,\infty}$)  so that NLS is still local and global well posed in suitable function spaces?
\end{que}
 The above question is also standing for the other dispersive equations. In \cite{KrSch14}, Krieger and Schlag considered a class of energy-supercritical NLW in 3D (higher dimensional cases are similar to 3D)
\begin{align}
 u_{tt}  -\Delta  u \pm  u^7 =0, \quad u(0,x)=u_0(x),\; u_t(0,x)=u_1(x)
\label{NLW}
\end{align}
and they obtained the following result (notice that $\dot H^{7/6}  \times \dot H^{1/6}$ is a critical space):

\begin{thm}
   Let $d=3$. There exist $(u_0, u_1) \in C^\infty \times C^\infty$ with
$$
\|(u_0, u_1)\|_{\dot H^{7/6}  \times \dot H^{1/6}}= \infty, \ \  \|(u_0, u_1)\|_{\dot H^{s} \times \dot H^{s-1}} < \infty, \  s>7/6
$$
such that there exists globally in
forward time as a $C^\infty$-smooth solution. These solutions are stable under a certain
class of perturbations.
\end{thm}

Moreover, Krieger and Schlag considered the global existence for some initial data out of $\dot H^{7/6}  \times \dot H^{1/6}$ and large enough in $\dot B^{7/6}_{2,\infty}  \times \dot B^{1/6}_{2,\infty}$  and   $L^\infty(|x|\geq 1)$ in the defocusing case.

 On the other hand, in the focusing case $\lambda>0$, it is known that NLS generates blow up solutions at finite time for a class of initial data with negative energy. Let us recall Glassey's result (cf. \cite{Gl77}):
Let $\lambda>0$, $\kappa > 2/d$ and $\kappa \in \mathbb{N}$, $u_0 \in \mathscr{S}$ satisfying
\begin{align}
E(u_0)  <0, \ \ \ \
    &  \Im \int x \overline{u}_0 \nabla u_0 <0,   \label{neg2}
\end{align}
then the solution of \eqref{NLS} blows up at some finite time $T_0>0$ in the sense
\begin{align}
\lim_{t\to T_0}\|\nabla u (t)\|_2 =\infty, \ \ \lim_{t\to T_0}\|u (t)\|_\infty=\infty.   \label{blowsup}
\end{align}
Recently Du, Wu and Zhang (cf. \cite{DuWuZh16}) obtained for the focusing NLS, if
$\kappa \geq 2/(d-2)$ and $\kappa \in \mathbb{N}$.  $u_0 \in H^s$, $s> s(\kappa)$ and
\begin{align}
E(u_0) \|u_0\|^2_2   <
     \left|  \Im \int   \overline{u}_0 \nabla u_0 \right|^2,   \label{neg3}
\end{align}
then there exists  $T_0 \leq \infty$  such that the solution of \eqref{NLS} blows up at $T_0$ in the sense that
\begin{align}
& T_0 <\infty, \ \ \lim_{t\to T_0}\| u (t)\|_{H^s} =\infty,   \label{blowsupneg}\\
& T_0=\infty, \ \ \lim_{t\to T_0}\|  u (t)\|_{p} =\infty, \ \ \forall p>2\kappa+2. \label{blowup3}
\end{align}
In view of Glassey's and Du, Wu and Zhang's results, it is impossible to set up the global well-posedness for the focusing NLS without any conditions on initial data. In fact, for the energy critical focusing NLS, Kenig and Merle \cite{KeMe06} obtained some optimal global well-posedness and scattering results in the radial case by using the profile decomposition together with the concentration compactness techniques.  A sharp condition for scattering of the radial 3D cubic focusing NLS was obtained by Holmer and Roudenko \cite{HoRo08}.

In this paper, we mainly study the global well posedness of
\eqref{NLS} with small initial data in $\alpha$-modulation spaces $M^{s,\alpha}_{2,1}$ and show the global well posedness of NLS in an $\alpha$-modulation spaces $M^{s,\alpha}_{2,1}$ and as an application of the well-posedness in $M^{s,\alpha}_{2,1}$,  we have a positive answer to Question \ref{question}.
$\alpha$-modulation spaces
$M_{p,q}^{s,\alpha}$ were introduced by Gr\"obner \cite{Gr92}, which can be regarded as  intermediate function spaces
to connect modulation and Besov spaces.  In  \cite{HaWa14}, we studied some standard properties including the dual space, embedding, scaling and algebraic structure of $\alpha$-modulation spaces and in current work we give their applications to NLS.
Now we recall the definition of $M_{2,1}^{s,\alpha}$.  Let $ Q_k^{\alpha}:= \langle k\rangle^{\frac{\alpha}{1-\alpha}}(k + [-C, C]^d) $ and
\begin{align} \label{alphamodulation}
 \|f\|_{M_{2,1}^{s,\alpha}}= \sum_{k\in \mathbb{Z}^d}  \langle
k\rangle^{\frac{s}{1-\alpha}}  \|  \widehat{f} \|_{L^2(Q_k^{\alpha} )}
\end{align}
In the case $\alpha=0$, $M^s_{2,1}:= M^{s,\alpha}_{2,1}$ is said to be a modulation space, cf. Feichtinger \cite{F}. $\alpha$-modulation spaces have been applied in pseudo-differential operators and related PDE in recent years, cf. \cite{BGOR07,BoNi062,BoNi06a,DaFor08,FHW11,Gr92,KoSuTo09,ST,SuWaZh15,T,ToWa11,WaHe07,WaHu07,WaZhGu06}.

Recall that the NLS has
the following equivalent form:
\begin{align}
& u(t)= S(t)u_0 +{\rm i}\mathscr{A} (|u|^{2\kappa}u), \label{NLSi} \\
& S(t)= e^{{\rm i} t\Delta}, \quad \mathscr{A}= \int^t_0 S(t-\tau) \cdot d\tau. \label{1.5a}
\end{align}
We have the following results.

\begin{thm} \label{Theorem1}
 Let $d\geq 1$, $\kappa \geq 2/d$ and $\kappa \in \mathbb{N}$, $s_\kappa = d\alpha/2 - \alpha / \kappa$ and $s>s_\kappa$.  Assume that $u_0 \in M^{s,\alpha}_{2,1}$. There exists $\varepsilon>0$ such that if $\|u_0\|_{M^{s,\alpha}_{2,1}} \leq \varepsilon$, then \eqref{NLSi} has a unique solution $u\in C(\mathbb{R}_+, M^{s,\alpha}_{2,1}) \cap X^{s,\alpha}_{\Delta}$\footnote{$X^{s,\alpha}_{\Delta}$ will be defined in Section
 \ref{Xsdelta}} and
 $$
\|u\|_{ C(\mathbb{R}_+, M^{s,\alpha}_{2,1}) \cap X^{s,\alpha}_{\Delta}} \lesssim \varepsilon.
 $$
 If $s<s_\kappa $,  then NLS is illposed in the sense that the solution map $u_0 \to u(t)$ is not $C^{2\kappa+1}$ in $M^{s,\alpha}_{2,1}$. Moreover, in the case $\alpha=0$, the solution map $u_0 \to u(t)$ is discontinuous at zero in $M^{s}_{2,1}$ if $s<s_k=0$.
\end{thm}

 We point out that the argument of Theorem \ref{Theorem1} implies the scattering operators of NLS carries a neighbourhood of zero in $M^{s,\alpha}_{2,1}$ into  $M^{s,\alpha}_{2,1}$ if the conditions of Theorem \ref{Theorem1} are satisfied.

Let us recall the definition of Besov spaces.
Let $s\in\mathbb{R}$, set $R_0:=\{\xi:|\xi|\leq 1\}, R_j:=\{\xi: 2^{j-1}\leq |\xi|<2^j\}$. Recall that the norm on Besov spaces $B^s_{2,q} \ (1\le q\le\infty)$ can be defined  as
\begin{equation*}
  \|f\|_{B^s_{2,q}}:=\left(\sum_{j\in\mathbb{Z}_+}2^{jsq}\|\widehat{f}\|_{L^2(R_j)}^q\right)^{1/q}<\infty
\end{equation*}
with usual modification for $q=\infty$.

\begin{thm}  \label{Corollary1}

Let $d\geq 1$, $\kappa > 2/d$ and $\kappa \in \mathbb{N}$,  $\alpha\in (0,1)$,  $d\alpha/2-\alpha/\kappa< s < s(\kappa)$, $M\gg 1$.  Then there exists $u_0 \in M^{s,\alpha}_{2,1}$ satisfying
\begin{align}
& \|u_0\|_2 \sim   M^{1/\kappa-d/2},   \ \       \|u_0\|_{M^{s,\alpha}_{2,1}} \geq  M^{(1-\alpha)/\kappa} ;    \label{largem} \\
&  \|u_0\|_\infty =\infty  , \ \ \|u_0\|_{B^{s(\kappa)}_{2,\infty}} =\infty \label{largeinfty}
\end{align}
such that \eqref{NLSi} has a unique solution $u\in C(\mathbb{R}, M^{s,\alpha}_{2,1}) \cap X^{s,\alpha}_{\Delta}$.
\end{thm}
Comparing the above result with Krieger and Schlag's Theorem 1.1, we see that the initial data $u_0 \in M^{s,\alpha}_{2,1}$ cannot be $C^\infty$ functions. However, noticing that $\|P_N u_0 -u_0 \|_{M^{s,\alpha}_{2,1}} \to 0$ as $N\to \infty$, where $P_N =\mathscr{F}^{-1} \chi_{|\xi|\leq N}\mathscr{F}$, we can consider $P_N u_0 \in C^\infty$ for $N\gg 1$ as the initial value so that the solution is $C^\infty$. We easily see that $P_N u_0 $  satisfies \eqref{largem} and by carefully choosing $u_0$, we have
\begin{align}
\|u_0\|_\infty \gg 1 , \ \ \|u_0\|_{B^{s(\kappa)}_{2,\infty}} \gg 1, \label{largeinftya}
\end{align}
see Section \ref{pfcor} for details. In the case $\alpha=0$, we can show that there exist  a class of initial data satisfying for any $M\gg 1$,
\begin{align}
& \|u_0\|_2 \sim   M^{1/\kappa-d/2},   \ \       \|u_0\|_{M^{s}_{2,1}} \geq  M^{1/\kappa} ;    \label{largem0} \\
&  M^{1/\kappa} \leq \|u_0\|_\infty <\infty, \ \ \|u_0\|_{B^{s(\kappa)}_{2,\infty}} =\infty  \label{largeinfty0}
\end{align}
such that \eqref{NLSi} has a unique solution $u\in C(\mathbb{R}, M^{s}_{2,1}) \cap X^{s,0}_{\Delta}$ for $0<s\ll 1$, $d\geq 1$, $\kappa>2/d$. The results in the case $\alpha =0 $ were essentially obtained in \cite{WaHe}.

   Theorem \ref{Corollary1} is particularly interesting to the focusing NLS with the energy super-critical power $\kappa>2/(d-2)$, which implies that there exist a class of initial data out of the critical Besov spaces $B^{s(\kappa)}_{2,\infty}$ so that NLS is globally well posed. Moreover, such a kind of data have neither small amplitudes nor small $\alpha$-modulation norm in  $M^{s,\alpha}_{2,1}$, $s=d\alpha/2-\alpha/\kappa +$.

In \cite{HaWa14a}, we announced a result to showed that NLS is globally well posed in $\alpha$-modulation spaces $M^{s,\alpha}_{2,1}$ if $s>\widetilde{s}_k$, where
$$
\widetilde{s}_k = \frac{d\alpha}{2} -  \frac{\alpha}{\kappa} + \frac{\alpha(1-\alpha)(n\kappa+2)}{2\kappa[(1+\alpha)\kappa+1-\alpha]}.
$$
$\widetilde{s}_k$ is not optimal if $\alpha\neq 0,1$.   Our Theorem \ref{Theorem1} has improved the result in \cite{HaWa14a} and except for the end point case $s=s_\kappa$, our result is sharp. The proof of Theorem \ref{Theorem1} relies upon $U_p$ and $V_p$ spaces first applied in \cite{HaHeKo09,KoTa05}, together with the  bilinear and Strichartz' estimates and variant $\alpha$-decompositions.

Throughout this paper,
$A\lesssim B$ stands for $A\le C B$, and $A\sim B$ denote $A\lesssim B$ and $B\lesssim A$, where $C$ is a positive constant which can be different at different places. We write $a\vee b=\max\,(a,b)$, $a\wedge b =\min\,(a,b)$. We denote by $s+= s+\varepsilon$, $0<\varepsilon \ll 1.$  For any $p\in [1,\infty]$, $p'$ will stand for
the dual number of $p$, i.e., $1/p+1/p'=1$.

The paper is organized as follows.
In Section \ref{functspace}, we introduce some function spaces and study their properties, such as duality, embedding, which are useful in the whole paper.
In Section \ref{bilinearest}, we get some bilinear estimates with respect to the $\alpha$-decompositions.
Global well-posedness for NLS in $\alpha$-modulation space is obtained in Sections \ref{multilin} and \ref{pfcor}.
In the last section, we show the ill-posedness for NLS in $\alpha$-modulation space.

\section{Function spaces} \label{functspace}

 Let $\mathscr{S}(\mathbb{R}^d)$  be the Schwartz space and $\mathscr{S}'(\mathbb{R}^d)$ be its  dual space.   Sobolev spaces
$
H^s(\mathbb{R}^d)=(I-\Delta)^{-s/2} L^2 , \ \   \dot H^s(\mathbb{R}^d)=( -\Delta)^{-s/2} L^2.
$
$L^q_{t} L^p_{x}(\mathbb{R}^{d+1})$ is the space-time Lebesgue space equipped with norm
\begin{equation*}
  \|f\|_{L^q_{t} L^p_{x}(\mathbb{R}^{d+1}) }:=\|\|f\|_{ L^p_x(\mathbb{R}^d)}\|_{L^q_t(\mathbb{R})}.
\end{equation*}

\subsection{Decompositions to frequency spaces}

Let
$\varphi$ be a smooth radial bump function supported in the ball $B(0,2)=\{\xi\in\mathbb{R}^d:|\xi|<2\}$
and $\varphi(\xi)=1$ in the unit ball $B(0,1)$. Put
\begin{equation}\label{varphi}
 \begin{cases}
 \varphi_0(\xi)=\varphi(\xi),
 \\
 \varphi_j(\xi)=\varphi(2^{-j}\xi)-\varphi(2^{-j+1}\xi), \quad j\in\mathbb{N},
 \end{cases}
\end{equation}
and
\begin{equation}\label{}
  \triangle_j:=\mathscr{F}^{-1}\varphi_j\mathscr{F}, \quad j\in\mathbb{N}\cup\{0\},
\end{equation}
which are said to be dyadic decomposition operators.  Put
\begin{equation}
\varphi_{k}^{\alpha}(\xi):=\varphi\left(\frac{\xi-\langle k\rangle^{\frac{\alpha}{1-\alpha}} k }{C\langle k \rangle^{\frac{\alpha}{1-\alpha}}}\right), \label{phialpha}
\end{equation}
and denote
\begin{equation*}
\eta_{ k }^{\alpha}(\xi):=\varphi_{ k }^{\alpha}(\xi)\left(\sum_{\boldsymbol{l}\in\mathbb{Z}^d}\varphi_{ l }^{\alpha}(\xi)\right)^{-1}.
\end{equation*}
We define an
operator sequence by
\begin{equation}
\square_{ k }^{\alpha}:=\mathscr{F}^{-1}\eta_{ k }^{\alpha}\mathscr{F}, \quad  k \in\mathbb{Z}^d,
\end{equation}
which are said to be the $\alpha$-decomposition operators. Formally, we have
$$
\sum_{j\in \mathbb{Z}_+} \triangle_j  = {\rm Id}, \ \ \    \sum_{k\in \mathbb{Z}^n} \Box^\alpha_k = {\rm Id}.
$$

The following are the embedding results between $\alpha$-modulation and Besov spaces.  One can refer to \cite{HaWa14}.
\begin{prop}[Embeddings]\label{embedding-mod-alpha}
There hold the following sharp embeddings.
\renewcommand{\labelenumi}{\roman{enumi}$)$}
\begin{enumerate}

\item[\rm (i)]
If $0\le\alpha<1,s_1\ge s_2$, then
$
  M^{s_1,\alpha}_{2,1}\subset B^{s_2}_{2,1}.
$

\item[\rm (ii)]
If $0\le\alpha<1,s_1\ge s_2+d(1-\alpha)/2$, then
$
  B^{s_1}_{2,1}\subset M^{s_2,\alpha}_{2,1}.
$
\end{enumerate}
\end{prop}

\subsection{$U^p$ and $V^p$ type spaces}
$U^p$ and $V^p$, as a development of Bourgain's spaces \cite{Bo93I,Bo93II} were first applied by Koch and Tataru in the study of NLS, cf. \cite{KoTa05,KoTa07,KoTa12}.    Let $\mathcal{Z}$ be the set of finite partitions $-\infty= t_0 <t_1<...< t_{K-1} < t_K =\infty$. Let $1\leq p <\infty$. For any $\{t_k\}^K_{k=0} \subset \mathcal{Z}$ and $\{\phi_k\}^{K-1}_{k=0} \in L^2$ with $\sum^{K-1}_{k=0} \|\phi_k\|^p_2=1$, $\phi_0=0$. A step function $a: \mathbb{R}\to L^2$ given by
$$
a= \sum^{K}_{k=1} \chi_{[t_{k-1}, t_k)} \phi_{k-1}
$$
is said to be a $U^p$-atom. All of the $U^p$ atoms is denoted by $\mathcal{A}(U^p)$.   The $U^p$ space is
$$
U^p=\left\{u= \sum^\infty_{j=1} c_j a_j : \ a_j \in \mathcal{A}(U^p), \ \ c_j \in \mathbb{C}, \ \ \sum^\infty_{j=1} |c_j|<\infty  \right\}
$$
for which the norm is given by
$$
\|u\|_{U^p}= \inf \left\{\sum^\infty_{j=1} |c_j| : \ \ u= \sum^\infty_{j=1} c_j a_j , \ \  \ a_j \in \mathcal{A}(U^p), \ \ c_j \in \mathbb{C} \right\}.
$$

We define $V^p$ as the normed space of all functions $v: \mathbb{R} \to L^2$ such that $\lim_{t\to \pm \infty} v(t)$ exist and for which the norm
$$
\|v\|_{V^p} := \sup_{\{t_k\}^K_{k=0} \in \mathcal{Z}} \left( \sum^K_{k=1} \|v(t_k)-v(t_{k-1})\|^p_{L^2}\right)^{1/p}
$$
is finite, where we use the convention that $v(-\infty) = \lim_{t\to -\infty} v(t)$ and $v(\infty)=0$ (here $v(\infty)$ and $\lim_{t\to  \infty} v(t)$ are different notations). Likewise, we denote by $V^p_-$ the subspace of all $v\in V^p$ so that $v(-\infty) =0$. Moreover, we define the closed subspace $V^p_{rc}$ $(V^p_{-,rc})$ as all of the right continuous functions in $V^p$ $(V^p_-)$.

We define
$$
U^p_{\Delta} := e^{\cdot \ {\rm i}\Delta} U^p, \ \ \|u\|_{U^p_{\Delta}} = \|e^{-{\rm i} t \Delta} u \|_{U^p}.
$$
$$
V^p_{\Delta} := e^{\cdot \ {\rm i}\Delta} V^p, \ \ \|u\|_{V^p_{\Delta}} = \|e^{-{\rm i} t \Delta} u \|_{V^p}.
$$
Similarly for the definition of $V^p_{rc, \Delta}$, $V^p_{-, \Delta}$, $V^p_{rc,-, \Delta}$.
We list some known results in $U^p$ and $V^p$ (cf. \cite{KoTa05,KoTa07,KoTa12,HaHeKo09}).

\begin{prop} \label{UVprop1}
{\rm (Embedding)} Let $1\leq p <q <\infty$. We have the following results.
\begin{itemize}
 \item[\rm (i)]  $U^p$ and $V^p$, $V^p_{rc}$, $V^p_{-}$, $V^p_{rc, -}$ are Banach spaces.

\item[\rm (ii)] $U^p\subset V^p_{rc, -} \subset U^q \subset L^\infty (\mathbb{R}, L^2)$. Every $u\in U^p$ is right continuous on $t\in \mathbb{R}$

\item[\rm (iii)] $V^p \subset V^q$,   $V^p_{-} \subset V^q_{-} $,   $V^p_{rc} \subset V^q_{rc} $,  $V^p_{rc, -} \subset V^q_{rc, -} $.

\item[\rm (iv)] $\dot X^{0, 1/2, 1} \subset U^2_{\Delta} \subset V^2_{rc,-,\Delta} \subset \dot X^{0, 1/2, \infty}$, where Besov type Bourgain's spaces $\dot X^{s, b, q}$ are defined by
$$
\|u\|_{\dot X^{s,b,q}} := \left\| \|\chi_{|\tau+\xi^2|\in [2^{j-1}, 2^j)} |\xi|^{s} |\tau+\xi^2|^{b} \widehat{u}(\tau,\xi) \|_{L^2_{\xi,\tau}}  \right\|_{\ell^q_{j\in \mathbb{Z}}}.
$$
\end{itemize}
\end{prop}

The following transference principle due to Hadac, Herr and Koch \cite{HaHeKo09} will be frequently used in the later.
\begin{prop}[{Transference principle}] \label{transf}
Suppose that the linear operator $T: L^2(\mathbb{R}^d) \times...\times L^2 (\mathbb{R}^d) \to L^q_tL^r_x(\mathbb{R}^{d+1})$ satisfies for some $1\le q,r\le\infty$,
\begin{equation}\label{transference-2222}
  \|T(({\rm e}^{{\rm i}t\Delta}\phi^{(\ell)})_{\ell=1}^n)\|_{L^q_tL^r_x(\mathbb{R}^{d+1})}\le C_1\prod_{\ell=1}^n\|\phi^{(\ell)}\|_{L^2},
\end{equation}
then we have
\begin{equation}\label{transference-qqqq}
  \|T((u^{(\ell)})_{\ell=1}^n)\|_{L^q_tL^r_x(\mathbb{R}^{d+1})}\le C_1\prod_{\ell=1}^n\|u^{(\ell)}\|_{U^q_{\Delta}}.
\end{equation}
Furthermore, if  for some $p>q$, there holds
\begin{equation}\label{transference-pppp}
  \|T((u^{(\ell)})_{\ell=1}^n)\|_{L^q_tL^r_x(\mathbb{R}^{d+1})}\le C_2\prod_{\ell=1}^n\|u^{(\ell)}\|_{U^p_{\Delta}},
\end{equation}
then we conclude that
\begin{equation}\label{transference-VVVVqqqq}
  \|T((u^{(\ell)})_{\ell=1}^n)\|_{L^q_tL^r_x(\mathbb{R}^{d+1})} \lesssim C_1\left(\ln\frac{C_2}{C_1}+1\right)\prod_{\ell=1}^n\|u^{(\ell)}\|_{V^q_{\Delta}}.
\end{equation}
\end{prop}

\subsection{ $X^{s,\alpha}_{\Delta}$ and $Y^{s,\alpha}_{\Delta}$} \label{Xsdelta}

In this subsection we introduce the working spaces $X^{s,\alpha}_{\Delta}$ and $Y^{s,\alpha}_{\Delta}$ and establish their duality. Moreover, we will give some embedding results between $X^{s_1,\alpha_1}_{\Delta}$ and $X^{s_2,\alpha_2}_{\Delta}$ ($Y^{s_1,\alpha_1}_{\Delta}$ and $Y^{s_2,\alpha_2}_{\Delta}$), which are necessary for us to make nonlinear estimates.

\begin{defin}\label{def-Z-Y}
\renewcommand{\labelenumi}{\roman{enumi}$)$}

Let $0\leq \alpha <1$, $s\in \mathbb{R}$. We define the following spaces:

\begin{align}
 X^{s,\alpha}  & =\{u\in\mathscr{S}'(\mathbb{R}^{d+1}): \  \|u\|_{X^{s,\alpha} }:=\sum_{ k \in\mathbb{Z}^d}\langle k \rangle^{\frac{s}{1-\alpha}}\|\square_{ k }^{\alpha}u\|_{U ^2} <\infty \},   \\
 Y^{s,\alpha}   & =\{u\in\mathscr{S}'(\mathbb{R}^{d+1}):    \|v\|_{Y^{s,\alpha}}:=\sup_{ k \in\mathbb{Z}^d}\langle k \rangle^{\frac{s}{1-\alpha}}\|\square_{ k }^{\alpha}v\|_{V ^2} <\infty\},\\
  X^{s,1}   & =\{u\in\mathscr{S}'(\mathbb{R}^{d+1}):
  \|u\|_{X^{s,1} }:=\sum_{j=0}^{\infty}2^{js}\|\triangle_ju\|_{U ^2} < \infty\},\\
Y^{s,1}   & =\{u\in\mathscr{S}'(\mathbb{R}^{d+1}):  \|v\|_{Y^{s,1} }:=\sup_{j\in\mathbb{Z}_+}2^{js}\|\triangle_jv\|_{V ^2}  < \infty\},\\
 \|u\|_{X^{s,\beta}_{\Delta}} & = \|S(-t) u\|_{X^{s,\beta}}, \ \  \|u\|_{Y^{s,\beta}_{\Delta}} = \|S(-t) u\|_{Y^{s,\beta}}, \ \ \beta \in [0,1].
\end{align}
For any time interval $I\subset \mathbb{R}$, we denote
\begin{align} \label{localXs}
\|u\|_{X^{s,\alpha}_{\Delta}  (I)} =\inf \{\|\widetilde{u}\|_{X^{s,\alpha}_{\Delta} } : \ \widetilde{u} \in X^{s,\alpha}_{\Delta} ,   \  \widetilde{u}(t) = u(t), \ \ \forall \ t \in I   \}.
\end{align}

\end{defin}

\begin{prop} \label{UVprop3}
{\rm (Duality)} Let $1\leq p   <\infty$, $1/p+1/p'=1$.  Then $(U^p)^* = V^{p'}$ in the sense that
\begin{align}
T: V^{p'}  \to (U^p)^*; \ \ T(v)=B(\cdot,v), \label{dual}
\end{align}
is an isometric mapping.  The bilinear form $B: U^p\times V^{p'}$ is defined in the following way: For a partition $\mathrm{t}:= \{t_k\}^K_{k=0} \in \mathcal{Z}$, we define
 \begin{align}
B_{\mathrm{t}} (u,v) = \sum^K_{k=1} \langle u(t_{k-1}), \ v(t_k)-v(t_{k-1})\rangle. \label{dual2}
\end{align}
Here $\langle \cdot, \cdot \rangle$ denotes the inner product on $L^2$. For  any $u\in U^p$, $v\in V^{p'}$, there exists a unique number $B(u,v)$ satisfying the following property. For any $\varepsilon>0$, there exists a partition $\mathrm{t}$ such that
$$
|B(u,v)- B_{\mathrm{t}'} (u,v)| <\varepsilon, \ \ \forall \  \mathrm{t}'    \supset \mathrm{t}.
$$
Moreover,
$$
|B(u,v)| \leq \|u\|_{U^p} \|v\|_{V^{p'}}.
$$
In particular, let $u\in V^1_{-}$ be absolutely continuous on compact interval, then for any $v\in V^{p'}$,
$$
 B(u,v) =\int \langle u'(t), v(t)\rangle dt.
$$
\end{prop}
The duality of $U^p$ and $V^{p'}$ obtained by Hadac, Herr and Koch \cite{HaHeKo09} is of importance for us to make nonlinear estimates for the dispersive equations. We need further consider its localized version with $\alpha$-decomposition, namely, the duality of $X^{s,\alpha}$ and $Y^{-s,\alpha}$.

\subsection{Duality of $X^{s,\alpha} $ and $Y^{-s, \alpha}$}

\begin{prop}
[{\rm Duality}]\label{duality-z-y}
Let $1\le p,q<\infty$, $\alpha\in [0,1)$. Then we have
\begin{gather}
  (X^{s,\alpha})^*=Y^{-s,\alpha}; \ \ \
  (X^{s,1})^*=Y^{-s,1}. \label{duality-z-y-1}
\end{gather}
  in the sense that
\begin{align}
T: Y^{-s, \alpha}    \to (X^{s,\alpha})^* ; \ \ T(v)=B(\cdot,v), \label{dualprop4}
\end{align}
is an isometric mapping, where the bilinear form $B(\cdot,\cdot)$ is defined in  Proposition \ref{UVprop3}. Moreover, we have
$$
|B(u,v)| \leq   \|u\|_{X^{s,\alpha}} \|v\|_{Y^{-s,\alpha}}.
$$
\end{prop}
{\bf Proof.} By the almost orthogonality, we see that
$$
B_{\mathrm{t}}(\Box^{\alpha}_k u, \Box^{\alpha}_l v) =0, \ \ |k- l| \geq C,
$$
For any $v\in  Y^{-s,\alpha}$, by Proposition \ref{UVprop3} and H\"older's inequality,  we have
\begin{align}
|B(u,v)| & = \left|\sum_{|l|\leq C }\sum_{k\in \mathbb{Z}^d} B (\Box^{\alpha}_k u, \Box^{\alpha}_{k+\ell} v)\right| \nonumber\\
 & \leq  \sum_{k\in \mathbb{Z}^d} \|\Box^{\alpha}_k u\|_{U^2} \|\Box^{\alpha}_{k+l} v\|_{V^2}  \lesssim  \|u\|_{X^{s, \alpha}} \|v\|_{Y^{-s,\alpha}} \nonumber
\end{align}
It follows that $Y^{-s, \alpha} \subset (X^{s,\alpha})^*$ and $\|v\|_{(X^{s,\alpha})^*} \lesssim  \|v\|_{Y^{-s,\alpha}}$.

Conversely,  considering the map
$$
X^{s,\alpha} \ni f \to \{\Box^{\alpha}_k f\} \in \ell^s_1 (\mathbb{Z}^d; U^2),
$$
where
$$
\ell^s_q (\mathbb{Z}^d; U^2):= \left\{\{f_k\}_{k\in \mathbb{Z}^d}: \ \  \|\{f_k\}\|_{\ell^s_q (\mathbb{Z}^d; U^2)} := \left\| \{\langle k\rangle^{s/(1-\alpha)}\|f_k\|_{U^2} \} \right\|_{\ell^q} <\infty \right\},
$$
we see that it is an isometric mapping from $X^{s,\alpha}$ into a subspace of $\ell^s_1 (\mathbb{Z}^d; U^2)$.  So,  $v \in (X^{s,\alpha})^*$ can be regarded as a continuous functional in a subspace of $\ell^s_1 (\mathbb{Z}^d; U^2)$. In view of Hahn-Banach Theorem, it can be extended onto $\ell^s_1 (\mathbb{Z}^d; U^2)$ (the extension is written as $\tilde{v}$) and its norm will be preserved.  In view of the well-known duality $(\ell^s_1 (\mathbb{Z}; X))^* = \ell^{-s}_{\infty} (\mathbb{Z}^d; X^*)$, we have
$$
(\ell^s_1 (\mathbb{Z}^d; U^2))^* = \ell^{-s}_{\infty} (\mathbb{Z}^d; V^2),
$$
and there exists $\{v_k\}_{k\in \mathbb{Z}} \in \ell^{-s}_{\infty} (\mathbb{Z}^d; V^2)$ such that
$$
\langle \tilde{v}, \ \{f_k\}\rangle = \sum_{k\in \mathbb{Z}^d} B(f_k, v_k), \ \ \forall \ \{f_k\} \in \ell^s_1 (\mathbb{Z}^d; U^2).
$$
Moreover, $\|v\|_{(X^{s,\alpha})^*}= \|\{v_k\}\|_{\ell^{-s,\alpha} (\mathbb{Z}^d; V^2)}.$  Hence,  for any $u\in X^{s,\alpha}$,
$$
\langle v, u\rangle = \langle \tilde{v}, \{\Box^{\alpha}_k u\} \rangle =  \sum_{k\in \mathbb{Z}^d} B(\Box^{\alpha}_k u, v_k).
$$
From $B_{\mathrm{t}}(\Box^{\alpha}_k u, \ v) = B_{\mathrm{t}}( u, \ \Box^{\alpha}_k v)$ we see that $B (\Box^{\alpha}_k u, \ v) = B ( u, \ \Box^{\alpha}_k v)$. It follows that
$$
v= \sum_{k\in \mathbb{Z}^d} \Box^{\alpha}_k v_k.
$$
Obviously, we have
$$
\|v\|_{Y^{-s,\alpha}} \leq \|\{v_k\}\|_{\ell^{-s,\alpha} (\mathbb{Z}; V^2)} = \|v\|_{(X^{s,\alpha})^*}.
$$
This proves $(X^{s,\alpha})^* \subset Y^{-s,\alpha}$.  $\hfill\Box$

Now we apply the duality to the norm calculation to the inhomogeneous part of the solution of NLS in $X^{s,\alpha}_{\Delta}$.
By Propositions \ref{UVprop3} and \ref{duality-z-y}, we see that
\begin{align}
\|\mathscr{A} (f) \|_{X^{s,\alpha}_{\Delta}} &  = \sup \left\{ \left|B\left(\int^t_0 e^{-{\rm i}s\Delta} f(s) ds, v \right) \right| : \  \|v\|_{Y^{-s,\alpha}_{p'}} \leq 1 \right\} \nonumber\\
& \leq \sup_{\|v\|_{Y^{-s,\alpha}} \leq 1} \left|\int_{  } \ \langle f (s), \  e^{{\rm i}s \Delta} v(s) \rangle ds  \right| \nonumber\\
& \leq \sup_{\|v\|_{Y^{-s,\alpha}_{ \Delta}} \leq 1} \left|\int_{  }\  \langle f (s), \  v(s) \rangle ds  \right|.  \label{normest}
\end{align}

\begin{cor}
 Let $0\le\alpha<1$. We have
\begin{align}\label{inhomogeneous-f-aaaa}
 &  \|\mathscr{A}f\|_{X^{s,\alpha}_{\Delta}}=\sup_{\|v\|_{Y^{-s,\alpha}_{\Delta}}\le 1}\left|\int_{\mathbb{R}}\langle f(t),v(t)\rangle dt\right|,
\\
 &  \|\mathscr{A}f\|_{X^{s,1}_{\Delta}}=\sup_{\|v\|_{Y^{-s,1}_{\Delta}}\le1}\left|\int_{\mathbb{R}}\langle f(t),v(t)\rangle dt\right|.  \label{inhomogeneous-f-1111}
\end{align}

\end{cor}

For the purpose of later use, we need:
\begin{thm}[{Embeddings}]\label{embedding-z-alpha}
The following embeddings hold true.
\renewcommand{\labelenumi}{\roman{enumi}$)$}
\begin{enumerate}
\item[\rm (i)]
If $0\le\alpha_1,\alpha_2<1,s_1\ge s_2+0\vee d(\alpha_1-\alpha_2)/2$, then
$
  X^{s_1,\alpha_1}\subset X^{s_2,\alpha_2}.
$

\item[\rm (ii)]
If $0\le\alpha<1,s_1\ge s_2+d(1-\alpha)/2$, then
$
  X^{s_1, 1}\subset X^{s_2,\alpha}.
$

\item[\rm (iii)]  If $s_1\ge s_2+0\vee d(\alpha_2-\alpha_1)/2$, then  $Y^{s_1,\alpha_1}\subset Y^{s_2,\alpha_2}$.

\item[\rm (iii)]  If $0\le\alpha<1$,  $s_1\ge s_2+  d(1-\alpha_1)/2$, then  $Y^{s_1,\alpha_1}\subset Y^{s_2,1}$.
\end{enumerate}
\end{thm}
\begin{proof}
For any $(l,\alpha_2)\in\mathbb{Z}^d\times[0,1)$, we define
\begin{equation*}
  \Lambda^{\alpha_1}[(l,\alpha_2)]
  :=\{k\in\mathbb{Z}^d:\square_k^{\alpha_1}\square_l^{\alpha_2}f\neq0,\;\forall f\in\mathscr{S}'(\mathbb{R}^d)\}.
\end{equation*}
One can check that
\begin{gather*}
  \#\Lambda^{\alpha_1}[(l,\alpha_2)]
  \lesssim
  1\vee\langle l\rangle^{\frac{d(\alpha_2-\alpha_1)}{1-\alpha_2}}; \ \ and \ \
   \langle k\rangle\sim\langle l\rangle^{\frac{1-\alpha_1}{1-\alpha_2}}
\end{gather*}
provided that $k\in\Lambda^{\alpha_1}[(l,\alpha_2)]$. Let $\{t_i\}_{i=0}^K\in\mathcal{Z}$, Plancherel's identity yields that
\begin{equation}\label{embedding-Y-1111}
  \begin{split}
     &
     \left(\sum_{i=1}^K\|\square_l^{\alpha_2}[u(t_i)-u(t_{i-1})]\|_2^2\right)^{\frac12}
      \\
      & \quad\le
     \left(\sum_{k\in\Lambda^{\alpha_1}[(l,\alpha_2)]}\sum_{i=1}^K\|\square_k^{\alpha_1}\square_l^{\alpha_2}[u(t_i)-u(t_{i-1})]\|_2^2\right)^{\frac12}
     \\
     &\quad\lesssim
     \langle k\rangle^{\frac{0\vee d(\alpha_2-\alpha_1)}{1-\alpha_1}}\sup_{\{t_i\}_{i=0}^K\in\mathcal{Z}}\left(\sum_{i=1}^K\|\square_k^{\alpha_1}[u(t_i)-u(t_{i-1})]\|_2^2\right)^{\frac12},
  \end{split}
\end{equation}
where the $k$ on the right-hand side of \eqref{embedding-Y-1111} belongs to $\Lambda^{\alpha_1}[(l,\alpha_2)]$. Therefore,
\begin{equation*}
  \begin{split}
     &
     \sup_{l\in\mathbb{Z}^d}\langle l\rangle^{\frac{s_2}{1-\alpha_2}}\sup_{\{t_i\}_{i=0}^K\in\mathcal{Z}}\left(\sum_{i=1}^K\|\square_l^{\alpha_2}[u(t_i)-u(t_{i-1})]\|_2^2\right)^{\frac12}
      \\
      &\quad\lesssim
     \sup_{k\in\mathbb{Z}^d}\langle k\rangle^{\frac{1}{1-\alpha_1}\big[s_2+\frac{d}{2}(\alpha_2-\alpha_1)\big]}\sup_{\{t_i\}_{i=0}^K\in\mathcal{Z}}\left(\sum_{i=1}^K\|\square_k^{\alpha_1}[u(t_i)-u(t_{i-1})]\|_2^2\right)^{\frac12},
  \end{split}
\end{equation*}
which implies the result of (iii).
As a consequence of duality, we can obtain the result as desired. For instance, we prove (i). If $\alpha_1\geq \alpha_2$ and $s_1 \geq s_2+ d(\alpha_1-\alpha_2)/2$, then we have $Y^{-s_2, \alpha_2} \subset Y^{-s_1,\alpha_1}$. It follows from the duality that
\begin{align}
\|u\|_{X^{s_2,\alpha_2}} & =\sup_{v\in (X^{s_2,\alpha_2})^*\setminus \{0\}} \frac{B(u,v)}{\|v\|_{(X^{s_2,\alpha_2})^*} } \nonumber\\
& = \sup_{v\in Y^{-s_2,\alpha_2} \setminus \{0\}} \frac{B(u,v)}{\|v\|_{Y^{-s_2,\alpha_2} } } \nonumber\\
 & \leq  C \sup_{v\in Y^{-s_1,\alpha_1} \setminus \{0\}} \frac{B(u,v)}{\|v\|_{Y^{-s_1,\alpha_1} } }=C\|u\|_{X^{s_1,\alpha_1}}.
\end{align}
 If $\alpha_1 <\alpha_2$ and $s_1 \geq s_2$, we can prove (i) in a similar way.
\end{proof}

\section{Bilinear Estimates} \label{bilinearest}

Using standard dual arguments, one has the following Strichartz estimates, which is useful to obtain the bilinear estimates, cf. \cite{HuCh14}.

\begin{prop}[{Strichartz}] $(q,r)$ is said to be a admissible pair if $q,r\ge 2,(q,r,d)\neq(2,\infty,2)$  and
\begin{equation}\label{condition-sigma-adimissible}
  \frac12-\frac1r-\frac{2}{dq}\ge 0.
\end{equation}
Let  $(q,r)$ and $(\widetilde{q},\widetilde{r})$ be two admissible pairs, then we have
\begin{subequations}
\begin{align}
\|\square_k^{\alpha}{\rm e}^{{\rm i}t\Delta}u_0\|_{L^q_tL^r_x}\lesssim\langle k\rangle^{\frac{d\alpha}{1-\alpha}\big(\frac12-\frac1r-\frac{2}{dq}\big)}\|\square_k^{\alpha}u_0\|_{L^2_x}.
\label{strichartz-111-schrodinger}
\end{align}
\end{subequations}
\end{prop}
Applying the transference principle, we immediately have
\begin{align}
\|\square_k^{\alpha}u \|_{L^q_tL^r_x}\lesssim\langle k\rangle^{\frac{d\alpha}{1-\alpha}\big(\frac12-\frac1r-\frac{2}{dq}\big)}\|\square_k^{\alpha}u\|_{U^q_\Delta}.
\label{strichartz-Uq-schrodinger}
\end{align}

\begin{lem} [\rm  General bilinear estimates] \label{Bilinear}

Let $D_1$ and $D_2$ are compact subsets of $\mathbb{R}^d$. Assume that ${\rm supp} \widehat{\varphi}_i \subset D_i$, $i=1,2$ and
\begin{align}
\lambda= \inf \{|\xi^{(1)}_1 -\xi^{(2)}_1 |: \ \xi^{(i)} \in D_i, \xi^{(i)}= (\xi^{(i)}_1,..., \xi^{(i)}_d), \ i=1,2 \}>0  \label{dist}
\end{align}
Denote $\overline{D} =\{\bar{\xi}: \ \exists \xi_1, (\xi_1, \bar{\xi}) \in D\} $.  Then we have
\begin{align}
 \|e^{{\rm i}t\Delta} \varphi_1 e^{{\rm i}t\Delta} \varphi_2 \|_{L^2_{x,t}} \lesssim \lambda^{-1/2}   (|\overline{D}_1| \wedge |\overline{D}_2|)^{1/2}  \|\varphi_1  \|_{2} \|\varphi_2  \|_{2}.   \label{bilinear}
\end{align}
\end{lem}

Lemma \ref{Bilinear} is known in the literatures, see for instance, \cite{Gr05} in 1D and \cite{Bo98,OzTs98} in higher dimensions in the dyadic version. However, the current version shows that the decay $\lambda^{-1/2}$ in  the bilinear estimate only depends on the distance of the supports of $\widehat{\varphi}_1$ and $\widehat{\varphi}_2$ in one direction.

{\bf Proof.} We have
\begin{align}
 \|e^{{\rm i}t\Delta} \varphi_1 e^{{\rm i}t\Delta} \varphi_2 \|_{L^2_{x,t}}  = \sup_{\|g\|_{L^2_{x,t}}=1} \int_{\mathbb{R}^{d+1}}  e^{{\rm i}t\Delta} \varphi_1 e^{{\rm i}t\Delta} \varphi_2 \ \overline{g(x,t)} dxdt.   \label{bilinear1}
\end{align}
By Plancherel's identity,
\begin{align}
\int_{\mathbb{R}^{d+1}}  e^{{\rm i}t\Delta} \varphi_1 e^{{\rm i}t\Delta} \varphi_2 \ \overline{g(x,t)} dxdt
&  =  \int_{\mathbb{R}^{d+1}}  (e^{-{\rm i}t|\xi_1|^2} \widehat{\varphi}_1) * (e^{-{\rm i}t|\xi_2|^2} \widehat{\varphi}_2) \ \overline{\widehat{g}(\xi,t)} d\xi dt   \nonumber \\
&  =  \int_{\mathbb{R}^{2d+1}}   e^{-{\rm i}t|\xi_1|^2} \widehat{\varphi}_1(\xi_1)  e^{-{\rm i}t|\xi_2|^2} \widehat{\varphi}_2(\xi_2) \ \overline{\widehat{g}(\xi_1+\xi_2,t)} d\xi_1d\xi_2 dt \nonumber \\
&  =  \int_{\mathbb{R}^{2d}}    \widehat{\varphi}_1(\xi_1)   \widehat{\varphi}_2(\xi_2) \ \overline{\widehat{g}(\xi_1+\xi_2, -|\xi_1|^2-|\xi_2|^2)} d\xi_1d\xi_2  \nonumber \\
&  \leq      \|\varphi\|_2    \|\varphi_2\|_2  \|\chi_{D_1}(\xi_1) \chi_{D_2}(\xi_2)\widehat{g}(\xi_1+\xi_2, -|\xi_1|^2-|\xi_2|^2)\|_{L^2_{\xi_1,\xi_2}}.
\label{bilinear2}
\end{align}
By Change of variables $\xi=\xi_1+\xi_2, \ \tau=-|\xi_1|^2-|\xi_2|^2$, we have
\begin{align}
  \|\chi_{D_1}(\xi_1) \chi_{D_2}(\xi_2)\widehat{g}(\xi_1+\xi_2, -|\xi_1|^2-|\xi_2|^2)\|_{L^2_{\xi_1,\xi_2}} \leq \lambda^{-1} | \overline{D}_2|
\label{bilinear3}
\end{align}
Since $\overline{D}_1$ has the same role as $\overline{D}_2$, we have the result, as desired. $\hfill\Box$\\

Since we also need to handle the bilinear version $u\overline{v}$, similar to Lemma \ref{Bilinear}, we have

\begin{lem} \label{Bilinear2}

  Assume that the conditions of Lemma \ref{Bilinear} are satisfied.  Then we have
\begin{align}
 \|e^{{\rm i}t\Delta} \varphi_1 \overline{e^{{\rm i}t\Delta} \varphi_2} \|_{L^2_{x,t}} \lesssim \lambda^{-1/2}   (|\overline{D}_1| \wedge |\overline{D}_2|)^{1/2}  \|\varphi_1  \|_{2} \|\varphi_2  \|_{2}.   \label{bilinear2-2}
\end{align}
\end{lem}
{\bf Proof.} In the same way as in the proof of Lemma \ref{Bilinear}, noticing that
\begin{align}
\int_{\mathbb{R}^{d+1}}  e^{{\rm i}t\Delta} \varphi_1 \overline{e^{{\rm i}t\Delta} \varphi_2} \ \overline{g(x,t)} dxdt
&  =  \int_{\mathbb{R}^{2d}}    \widehat{\varphi}_1(\xi_1)   \widehat{\varphi}_2(\xi_2) \ \overline{\widehat{g}(\xi_1 -\xi_2,  |\xi_1|^2-|\xi_2|^2)} d\xi_1d\xi_2, \label{bilinear2-3}
\end{align}
we can repeat the argument as in the proof of Lemma \ref{Bilinear} to have the result, as desired. $\hfill\Box$

\begin{cor}[$\Box^\alpha_k$-Bilinear estimates I]\label{bilinear-alpha-schrodinger}
Let $0\le\alpha<1, \ l\in\mathbb{Z}^d, \ |l|\gg1$. $|k| \leq |l|$. Suppose that
\begin{align}
|\langle k  \rangle^{\alpha/(1-\alpha)} k  - \langle l \rangle^{\alpha/(1-\alpha)} l| \gtrsim  \langle l \rangle^{1/(1-\alpha)}. \label{constraint1}
\end{align}
Then we have
\begin{gather}
  \|\square_{k}^{\alpha}u\square_{l }^{\alpha}v\|_{L^2_{t,x}}\lesssim
  {\langle k \rangle}^{ (d-1)\alpha/2(1-\alpha)} {\langle l \rangle }^{-1/2(1-\alpha)} \ln  \langle l\rangle
  \|\square_{k }^{\alpha}u\|_{V^2_{\Delta}}
  \|\square_{l }^{\alpha}v\|_{V^2_{\Delta}}; \label{bilinear-u-v}
  \\
  \|\square_{k }^{\alpha}\overline{u}\square_{l }^{\alpha}\overline{v}\|_{L^2_{t,x}}\lesssim
  {\langle k \rangle}^{ (d-1)\alpha/2(1-\alpha)} {\langle l \rangle }^{-1 /2(1-\alpha)} \ln  \langle l\rangle
  \|\square_{-k }^{\alpha}u\|_{V^2_{\Delta}}
  \|\square_{-l }^{\alpha}v\|_{V^2_{\Delta}}; \label{bilinear-u-conj-v-conj}
  \\
  \|\overline{\square_{k }^{\alpha}u}\square_{l }^{\alpha}v\|_{L^2_{t,x}}\lesssim
  {\langle k \rangle}^{ (d-1)\alpha/2(1-\alpha)} {\langle l \rangle }^{-1 /2(1-\alpha)} \ln  \langle l\rangle
  \|\square_{k }^{\alpha}u\|_{V^2_{\Delta}}
  \|\square_{l }^{\alpha}v\|_{V^2_{\Delta}}; \label{bilinear-u-longconj-v}
  \\
  \|\square_{k }^{\alpha}u\overline{\square_{l }^{\alpha}v}\|_{L^2_{t,x}}\lesssim
  {\langle k \rangle}^{ (d-1)\alpha/2(1-\alpha)} {\langle l \rangle }^{-1 /2(1-\alpha)} \ln  \langle l\rangle
  \|\square_{k }^{\alpha}u\|_{V^2_{\Delta}}
  \|\square_{l }^{\alpha}v\|_{V^2_{\Delta}}. \label{bilinear-u-v-longconj}
\end{gather}
Moreover, if $2^j\ll \langle l\rangle^{1/(1-\alpha)}$, then we have
\begin{gather}
  \|\overline{\triangle_ju}\square_{l }^{\alpha}v\|_{L^2_{t,x}}\lesssim
  2^{j\frac{d-1}{2}}  {\langle l \rangle }^{-1/2(1-\alpha)} \ln  \langle l\rangle
  \|\triangle_ju\|_{V^2_{\Delta}}
  \|\square_{l }^{\alpha}v\|_{V^2_{\Delta}}. \label{bilinear-tri-u-v}
\end{gather}
\end{cor}

Notice that for $k \in\mathbb{Z}^d, \ |k|\ll |l|$, or  $l_1k_1 <0$ with $|l_1| =\max_{1\leq j \leq d} |l_j|$, \eqref{constraint1} holds.

{\bf Proof.} By H\"older's inequality and \eqref{strichartz-Uq-schrodinger},
\begin{align}
  \|\square_{k}^{\alpha}u\square_{l }^{\alpha}v\|_{L^2_{t,x}} & \leq   \|\square_{k}^{\alpha}u \|_{L^4_{t,x}} \|\square_{k}^{\alpha}v \|_{L^4_{t,x}} \nonumber\\
   & \lesssim \langle k\rangle^{(d\alpha/4-\alpha/2)/(1-\alpha)}  \langle l\rangle^{(d\alpha/4-\alpha/2)/(1-\alpha)} \|\square_{k}^{\alpha} u\|_{U^4_{\Delta}} \|\square_{k}^{\alpha} v\|_{U^4_{\Delta}}.
    \label{bilinear-u-vcor}
\end{align}
Using Lemma \ref{Bilinear}, we have
\begin{align}
  \|\square_{k}^{\alpha}u\square_{l }^{\alpha}v\|_{L^2_{t,x}}  \leq   {\langle k \rangle}^{ (d-1)\alpha/2(1-\alpha)} {\langle l \rangle }^{-1/2(1-\alpha)}   \|\square_{k}^{\alpha} u\|_{U^2_{\Delta}} \|\square_{k}^{\alpha} v\|_{U^2_{\Delta}}.
    \label{bilinear-u-vcor2}
\end{align}
By transference principle, we have \eqref{bilinear-u-v}.  Noticing that $|\square_{k }^{\alpha}\overline{u}\square_{l }^{\alpha}\overline{v}|= |\square_{-k }^{\alpha} {u}\square_{-l }^{\alpha} {v}|$, we immediately have \eqref{bilinear-u-conj-v-conj}.

\begin{cor}[$\Box^\alpha_k$-Bilinear estimate II]\label{bilinear-alpha-schrodinger-222}
Let $0\le\alpha<1, \ l\in\mathbb{Z}^d, \ |l |\gg 1$, $|k| \leq  |l|$. Suppose that
\begin{align}
|\langle k  \rangle^{\alpha/(1-\alpha)} k  + \langle l \rangle^{\alpha/(1-\alpha)} l| \gtrsim  \langle l \rangle^{1/(1-\alpha)}. \label{constraint2}
\end{align}

Then we have
\begin{gather}
  \|\square_{k }^{\alpha}\overline{u}\square_{l }^{\alpha}v\|_{L^2_{t,x}}\lesssim
  {\langle k \rangle}^{ (d-1)\alpha/2(1-\alpha)} {\langle l \rangle }^{-1 /2(1-\alpha)} \ln  \langle l\rangle
  \|\square_{-k }^{\alpha}u\|_{V^2_{\Delta}}
  \|\square_{l }^{\alpha}v\|_{V^2_{\Delta}}; \label{bilinear-u-conj-v}
  \\
  \|\square_{k }^{\alpha}u\square_{l }^{\alpha}\overline{v}\|_{L^2_{t,x}}\lesssim
  {\langle k \rangle}^{ (d-1)\alpha/2(1-\alpha)} {\langle l \rangle }^{-1 /2(1-\alpha)} \ln  \langle l\rangle
  \|\square_{k }^{\alpha}u\|_{V^2_{\Delta}}
  \|\square_{-l }^{\alpha}v\|_{V^2_{\Delta}}; \label{bilinear-u-v-conj}
  \\
  \|\overline{\square_{k }^{\alpha}u}\square_{l }^{\alpha}\overline{v}\|_{L^2_{t,x}}\lesssim
  {\langle k \rangle}^{ (d-1)\alpha/2(1-\alpha)} {\langle l \rangle }^{-1 /2(1-\alpha)} \ln  \langle l\rangle
  \|\square_{k }^{\alpha}u\|_{V^2_{\Delta}}
  \|\square_{-l }^{\alpha}v\|_{V^2_{\Delta}}; \label{bilinear-u-longconj-v-conj}
  \\
  \|\square_{k }^{\alpha}\overline{u}\overline{\square_{l }^{\alpha}v}\|_{L^2_{t,x}}\lesssim
   {\langle k \rangle}^{ (d-1)\alpha/2(1-\alpha)} {\langle l \rangle }^{-1 /2(1-\alpha)} \ln  \langle l\rangle
  \|\square_{-k }^{\alpha}u\|_{V^2_{\Delta}}
  \|\square_{l }^{\alpha}v\|_{V^2_{\Delta}}. \label{bilinear-u-conj-v-longconj}
\end{gather}
Moreover, if $2^j\ll \langle l\rangle^{1/(1-\alpha)}$, then we have
\begin{gather}
  \|\overline{\triangle_ju}\square_{l }^{\alpha}\overline{v}\|_{L^2_{t,x}}\lesssim
  2^{j\frac{d-1}{2}}  {\langle l \rangle }^{-1 /2(1-\alpha)} \ln  \langle l\rangle
  \|\triangle_ju\|_{V^2_{\Delta}}
  \|\square_{-l }^{\alpha}v\|_{V^2_{\Delta}}. \label{bilinear-tri-u-v-conj}
\end{gather}
\end{cor}

Noticing that if $ |k | \ll |l |$,  or $k_1l_1>0$ with  $|l_1| =\max_{1\leq j \leq d} |l_j|$,  \eqref{constraint2} holds.

\section{Multi-linear Estimates} \label{multilin}
Let us write
\begin{align}
  \mathscr{L} (u,v) =\int_{\mathbb{R}^{d+1}} |u|^{2\kappa} u \overline{v} dxdt .
\label{nonlinear1}
\end{align}

\begin{lem} \label{Nonlinear}
We have
 \begin{align}
|\mathscr{L} (u,v)| \lesssim  \prod^{2\kappa+1}_{j=1}  \|u\|^{2\kappa+1}_{X^{s,\alpha}_{\Delta}} \|v\|_{Y^{-s,\alpha}_{\Delta}}  \label{nonlinear0}
\end{align}

\end{lem}
For the proof of Lemma \ref{Nonlinear}, we will only consider the case $s=s_k+$, since the case of larger $s$ is easier to handle than that of $s=s_k$.
For convenience, we write
\begin{align}
 u_1=u_3=...=u_{2\kappa+1} =u, \ \ u_2=u_4=...=u_{2\kappa}= \overline{u}.
\label{nonlinear2}
\end{align}
Applying the $\alpha$-decompositions,
\begin{align}
  \mathscr{L} (u,v) = \sum_{  k^{(1)},...,k^{(2\kappa+1)} \in \mathbb{Z}^{d}} \int_{\mathbb{R}^{d+1}}  \prod^{2\kappa+1}_{\ell=1} \Box^\alpha_{k^{(\ell)}} u_{\ell}  \overline{v} dxdt .
\label{nonlinear3}
\end{align}
We further denote for given $k^{(1)},...,k^{(2\kappa+1)}$,
\begin{align}
K_j (k^{(1)},...,k^{(2\kappa+1)})= \left| \sum^{2\kappa+1}_{\ell=1} \langle k^{(\ell)} \rangle^{\alpha/(1-\alpha)} k^{(\ell)}_j \right| , \ \ j=1,...,2\kappa+1, \\ K:= K (k^{(1)},...,k^{(2\kappa+1)}) =\max_{j=1,...,d} K_j (k^{(1)},...,k^{(2\kappa+1)}).
\end{align}
Put
\begin{align}
\Omega_0 &  = \left\{(k^{(1)},..., k^{(2\kappa+1)}) \in \mathbb{Z}^{d(2\kappa+1)}: K \lesssim  \bigvee^{2\kappa+1}_{\ell=1} \langle k^{(\ell)} \rangle^{\alpha /(1-\alpha)} \right\},
   \label{kell1}\\
\Omega_1 & = \left\{(k^{(1)},..., k^{(2\kappa+1)}) \in \mathbb{Z}^{d(2\kappa+1)}: K \sim  \bigvee^{2\kappa+1}_{\ell=1} \langle k^{(\ell)} \rangle^{1/(1-\alpha)} \right\},
   \label{kell2}\\
\Omega_2 & = \left\{(k^{(1)},..., k^{(2\kappa+1)}) \in \mathbb{Z}^{d(2\kappa+1)}: \bigvee^{2\kappa+1}_{\ell=1} \langle k^{(\ell)} \rangle^{\alpha/(1-\alpha)}  \ll  K \ll  \bigvee^{2\kappa+1}_{\ell=1} \langle k^{(\ell)} \rangle^{1/(1-\alpha)} \right\}.
   \label{kell3}
\end{align}
We can divide $\mathscr{L} (u,v)$ into the following three parts:
\begin{align}
  \mathscr{L}_i (u,v) = \sum_{  (k^{(1)},...,k^{(2\kappa+1)}) \in \Omega_i } \int_{\mathbb{R}^{d+1}}  \prod^{2\kappa+1}_{\ell=1} \Box^\alpha_{k^{(\ell)}} u_{\ell}  \overline{v} dxdt .
\label{nonlinear3i}
\end{align}
For convenience, we will denote $(k^{(\ell)}):=(k^{(\ell)})^{2\kappa+1}_{\ell=1} = (k^{(1)},..., k^{(2\kappa +1)})$ in the sequel.

\subsection{Estimates of $\mathscr{L}_1 (u,v)$}

 Now we consider the estimate  of $\mathscr{L}_1 (u,v)$, which is easier than those of $\mathscr{L}_0 (u,v)$ and $\mathscr{L}_2 (u,v)$.  Using the $\alpha$-decomposition, we can rewrite $\mathscr{L}_1 (u,v)$ as
\begin{align}
  \mathscr{L}_1 (u,v) = \sum_{ (k^{(\ell)}) \in \Omega_1, \ k\in \mathbb{Z}^d} \int_{\mathbb{R}^{d+1}}  \prod^{2\kappa+1}_{\ell=1} \Box^\alpha_{k^{(\ell)}} u_{\ell}  \overline{\Box^\alpha_{k} v} dxdt .
\label{nonlinear31}
\end{align}
So, one needs to control the right hand side of  \eqref{nonlinear31} and we prove that
\begin{lem} \label{Nonlinearestim1}

\begin{align}
 |\mathscr{L}_{1} (u,v)| \lesssim   \|   v  \|_{Y^{-s, \alpha}_\Delta}   \| u \|^{2\kappa+1}_{X^{s,\alpha}_\Delta}.      \label{nonlinearestim1}
\end{align}
\end{lem}

. First, we have

\begin{lem} \label{Nonlinearnonzero}
If $\mathscr{L}_{1} (u,v)$ in \eqref{nonlinear31} does not equal to zero, we must have
\begin{align}
 \langle k \rangle^{\alpha/(1-\alpha)} (k_j -C) \leq  \sum^{2\kappa+1}_{\ell =1} \langle k^{(\ell)} \rangle^{\alpha/(1-\alpha)} (k^{(\ell)}_j + C),   \label{Nonlinearnonzero1}\\
 \langle k \rangle^{\alpha/(1-\alpha)} (k_j + C) \geq  \sum^{2\kappa+1}_{\ell =1} \langle k^{(\ell)} \rangle^{\alpha/(1-\alpha)} (k^{(\ell)}_j - C).   \label{Nonlinearnonzero2}
\end{align}
\end{lem}
{\bf Proof.} Using the fact that
$$
 \int_{\mathbb{R}^{d+1}}  \prod^{2\kappa+1}_{\ell=1} \Box^\alpha_{k^{(\ell)}} u_{\ell}  \overline{\Box^\alpha_{k} v} dxdt  = \int_{\mathbb{R}}  \mathscr{F} \left(\prod^{2\kappa+1}_{\ell=1} \Box^\alpha_{k^{(\ell)}} u_{\ell}  \overline{\Box^\alpha_{k} v}\right)(0,t)  dt,
$$
we easily see the result. $\hfill\Box$

According to the duality argument,  one needs to take $\ell^\infty$ norm on $k\in \mathbb{Z}^d$,  we must remove the summation on $k\in \mathbb{Z}^d$ in the right hand side of \eqref{nonlinear3}. For convenience, we write
\begin{align}
\Lambda((k^{(\ell)})^{2\kappa+1}_{\ell=1}) = \{k \in \mathbb{Z}^{d} : \eqref{Nonlinearnonzero1} \ {\rm and}  \ \eqref{Nonlinearnonzero2} \ {\rm are \ satisfied}\}. \label{lambda1}
\end{align}
So, we can furhter rewrite $\mathscr{L}_1 (u,v)$ as
\begin{align}
  \mathscr{L}_1 (u,v) = \sum_{ ((k^{(\ell)})^{2\kappa+1}_{\ell=1} \in \Omega_1, \ k\in \Lambda ((k^{(\ell)})^{2\kappa+1}_{\ell=1}) } \int_{\mathbb{R}^{d+1}}  \prod^{2\kappa+1}_{\ell=1} \Box^\alpha_{k^{(\ell)}} u_{\ell}  \overline{\Box^\alpha_{k} v} dxdt .
\label{nonlinear311}
\end{align}

For convenience, we denote by $\#A$ the number of the elements in $A$. It seems necessary to calculate $\# \Lambda ((k^{(\ell)})^{2\kappa+1}_{\ell=1})$.

\begin{lem} \label{Nonlinearnumber}
Suppose that $(k^{(\ell)})^{2\kappa+1}_{\ell=1} \in \Omega_1$,  $k\in \Lambda((k^{(\ell)})^{2\kappa+1}_{\ell=1})$. Then we have
\begin{align}
 K\sim  \langle k \rangle^{1/(1-\alpha)}, \ \ \# \Lambda ((k^{(\ell)})^{2\kappa+1}_{\ell=1}) \lesssim  1.
   \label{Nonlinearnumber2}
\end{align}

\end{lem}

{\bf Proof.} It is a straightforward consequence of \eqref{Nonlinearnonzero1} and \eqref{Nonlinearnonzero2}. $\hfill \Box$\\

{\bf Step 1.} We consider the estimates of $\mathscr{L}_1 (u,v)$. We can assume that $K\gg 1$, since all of the summations are finite terms if $K\lesssim 1$.  Since $k^{(1)}, k^{(3)},...,k^{(2\kappa+1)}$ have the equal positions, also  $k^{(2)}, k^{(4)},...,k^{(2\kappa)}$ play the same roles,  we can assume without loss of generality that
\begin{align}
  |k^{(1)}| \geq |k^{(3)}|\geq ...\geq |k^{(2\kappa+1)}|, \ \ |k^{(2)}| \geq |k^{(4)}|\geq ...\geq |k^{(2\kappa)}|. \label{Numberorder}
\end{align}
By \eqref{Nonlinearnonzero1} and \eqref{Nonlinearnonzero2}, for $(k^{(1)},...,k^{(2\kappa+1)}) \in \Omega_1 $, there exists $j\in \{1,...,d\}$ such that
\begin{align}
\langle k_j \rangle  \sim  \langle k \rangle  \sim  \langle k^{(1)} \rangle \sim  \langle k^{(1)}_j \rangle, \ or  \  \langle k_j \rangle  \sim \langle k \rangle  \sim  \langle k^{(2)} \rangle \sim \langle k^{(2)}_j \rangle.
\end{align}

{\it Case 1.} We consider the case
\begin{align}
\langle k_j \rangle  \sim  \langle k \rangle  \sim  \langle k^{(1)} \rangle \sim  \langle k^{(1)}_j \rangle. \label{nonlinear4}
\end{align}
We can assume that $j=1$ in \eqref{nonlinear4}.

{\it Case 1.1.} We consider the case $|k^{(2)}|, |k^{(3)}| \ll |k^{(1)}|$.  In $\mathscr{L}_1(u,v)$, using H\"older's inequality, bilinear and Strichartz' estimates, we can bound
\begin{align}
& \int_{\mathbb{R}^{d+1}}  \prod^{2\kappa+1}_{\ell=1} \Box^\alpha_{k^{(\ell)}} u_{\ell} \overline{\Box^\alpha_{k} v} dxdt \nonumber\\
& \leq \|\overline{\Box^\alpha_{k} v} \Box^\alpha_{k^{(3)}} u_{3} \|_{L^2_{x,t}} \|\Box^\alpha_{k^{(1)}} u_1  \Box^\alpha_{k^{(2)}} u_{2} \|_{L^2_{x,t}}  \prod^{2\kappa+1}_{\ell=4} \|\Box^\alpha_{k^{(\ell)}} u_{\ell}\|_{L^\infty_{x,t}} \nonumber\\
& \lesssim  \ln \langle k^{(1)}\rangle \ \langle k^{(1)} \rangle^{-1/(1-\alpha)}  \langle k^{(2)} \rangle^{(d-1)\alpha/2(1-\alpha)} \langle k^{(3)} \rangle^{(d-1)\alpha/2(1-\alpha)} \| \Box^\alpha_{k} v  \|_{V^2_\Delta} \|\Box^\alpha_{k^{(3)}} u  \|_{U^2_\Delta} \nonumber\\
& \quad \times \|\Box^\alpha_{k^{(1)}} u  \|_{U^2_\Delta} \| \Box^\alpha_{-k^{(2)}} u  \|_{U^2_\Delta} \prod^{2\kappa+1}_{\ell=4} \langle k^{(\ell)} \rangle^{d\alpha/2(1-\alpha)} \|\Box^\alpha_{\pm k^{(\ell)}} u_{\ell}\|_{U^2_\Delta } \nonumber\\
& \lesssim  \ln \langle k^{(1)}\rangle \ \langle k^{(1)} \rangle^{-1/(1-\alpha)}  \langle k^{(2)} \rangle^{((d-1)\alpha/2 -s)/(1-\alpha)} \langle k^{(3)} \rangle^{((d-1)\alpha/2 -s)/(1-\alpha)}  \langle  k  \rangle^{-s/(1-\alpha)}\| \Box^\alpha_{k} v  \|_{V^2_\Delta}   \nonumber\\
& \quad \times ( \langle k^{(4)} \rangle \cdot  \langle k^{(5)} \rangle )^{(\kappa-1)(d\alpha/2-s)/(1-\alpha)} \prod^{2\kappa+1}_{\ell=1} \langle k^{(\ell)} \rangle^{s/(1-\alpha)} \|\Box^\alpha_{\pm k^{(\ell)}} u_{\ell}\|_{U^2_\Delta },
\end{align}
where we used \eqref{Numberorder}. Further, by \eqref{Numberorder},
$$
 \langle k^{(4)} \rangle ^{(\kappa-1)(d\alpha/2-s)/(1-\alpha)}  \langle k^{(2)} \rangle^{((d-1)\alpha/2 -s)/(1-\alpha)} \leq   \langle k^{(2)} \rangle ^{(\kappa(d\alpha/2-s) -\alpha/2)/(1-\alpha)}.
$$
 Notice that
\begin{align}
  \kappa (d\alpha/2-s) -\alpha/2  \geq 0.
\end{align}
By \eqref{Numberlambda2}
we have $\# \Lambda(k^{(1)},...,k^{(2\kappa+1)}) \lesssim 1$ if $(k^{(1)},...,k^{(2\kappa+1)}) \in \Omega_1$.  If $|k^{(2)}|, |k^{(3)}| \ll |k^{(1)}|$ in the summation of $\mathscr{L}_1 (u,v)$, noticing that $-(1+\alpha) + 2\kappa(d\alpha/2-s)<0$,  we have
\begin{align}
\mathscr{L}_1 (u,v)
& \lesssim  \sum_{ ((k^{(\ell)})^{2\kappa+1}_{\ell=1} \in \Omega_1, \ k\in \Lambda ((k^{(\ell)})^{2\kappa+1}_{\ell=1}) } \ln \langle k^{(1)}\rangle \ \langle k^{(1)} \rangle^{-(1+\alpha)/(1-\alpha) + 2\kappa(d\alpha/2-s)/(1-\alpha)}   \nonumber\\
& \quad \times \langle k  \rangle^{- s/ (1-\alpha)} \| \Box^\alpha_{k} v  \|_{V^2_\Delta} \prod^{2\kappa+1}_{\ell=1} \langle k^{(\ell)} \rangle^{s/ (1-\alpha)} \|\Box^\alpha_{\pm k^{(\ell)}} u \|_{U^2_\Delta } \nonumber\\
& \lesssim   \|   v  \|_{Y^{-s}_\Delta}  \sum_{  k^{(1)},...,k^{(2\kappa+1)}\in \mathbb{Z}^d }
  \prod^{2\kappa+1}_{\ell=1} \langle k^{(\ell)} \rangle^{s/ (1-\alpha)} \|\Box^\alpha_{\pm k^{(\ell)}} u \|_{U^2_\Delta } \nonumber\\
& \lesssim   \|   v  \|_{Y^{-s,\alpha}_\Delta}   \| u \|^{2\kappa+1}_{X^{s,\alpha}_\Delta}.
\end{align}

{\it Case 1.2.} We consider the case $|k^{(2)}|\ll  |k^{(3)}| \sim |k^{(1)}|$ or $|k^{(3)}|\ll  |k^{(2)}| \sim |k^{(1)}|$. It suffices to consider the former case. Using H\"older's,  bilinear and Strichartz estimates,
\begin{align}
& \int_{\mathbb{R}^{d+1}}  \prod^{2\kappa+1}_{\ell=1} \Box^\alpha_{k^{(\ell)}} u_{\ell} \overline{\Box^\alpha_{k} v}  dxdt \nonumber\\
& \leq \|\overline{\Box^\alpha_{k} v} \Box^\alpha_{k^{(2)}} u_{2} \|_{L^2_{x,t}} \|\Box^\alpha_{k^{(1)}} u_1 \|_{L^4_{x,t}} \|\Box^\alpha_{k^{(3)}} u_{3} \|_{L^4_{x,t}}  \prod^{2\kappa+1}_{\ell=4} \|\Box^\alpha_{k^{(\ell)}} u_{\ell}\|_{L^\infty_{x,t}} \nonumber\\
& \lesssim  \ln \langle k^{(1)}\rangle \  \langle k^{(1)} \rangle^{-1/2(1-\alpha)}  \langle k^{(2)} \rangle^{(d-1)\alpha/2(1-\alpha)}  \langle k^{(1)} \rangle^{(d-2)\alpha/2(1-\alpha) } \| \Box^\alpha_{k} v  \|_{V^2_\Delta} \|\Box^\alpha_{k^{(3)}} u  \|_{U^2_\Delta} \nonumber\\
& \quad \times \|\Box^\alpha_{k^{(1)}} u  \|_{U^2_\Delta} \| \Box^\alpha_{-k^{(2)}} u  \|_{U^2_\Delta} \prod^{2\kappa+1}_{\ell=4} \langle k^{(\ell)} \rangle^{d\alpha/2(1-\alpha)} \|\Box^\alpha_{\pm k^{(\ell)}} u_{\ell}\|_{U^2_\Delta } \nonumber\\
& \lesssim  \ln \langle k^{(1)}\rangle \  \langle k^{(1)} \rangle^{(-1/2 -\alpha +\kappa(d\alpha/2-s) )/(1-\alpha)}  \langle k^{(2)} \rangle^{( -\alpha/2 +\kappa(d\alpha/2-s) )/(1-\alpha)}
\langle k \rangle^{-s/(1-\alpha)} \| \Box^\alpha_{k} v  \|_{V^2_\Delta}    \nonumber\\
& \quad \times \prod^{2\kappa+1}_{\ell=1} \langle k^{(\ell)} \rangle^{s/(1-\alpha)} \|\Box^\alpha_{\pm k^{(\ell)}} u_{\ell}\|_{U^2_\Delta }.
\end{align}
Noticing that $|k^{(1)}| \geq |k^{(2)}|$,
$$
-\alpha/2 +\kappa(d\alpha/2-s) >0, \ \ -1/2 -3\alpha/2 + 2\kappa(d\alpha/2-s) <0,
$$
and using a similar way to Case 1.1,  we can estimate $\mathscr{L}_1 (u,v)$ by
\begin{align}
\mathscr{L}_1 (u,v)
& \lesssim   \|   v  \|_{Y^{-s,\alpha}_\Delta}   \| u \|^{2\kappa+1}_{X^{s,\alpha}_\Delta}.
\end{align}

{\it Case 1.3.} We consider the case $|k^{(2)}|\sim  |k^{(3)}| \sim |k^{(1)}|$. Using H\"older's and Strichartz estimates,
\begin{align}
& \int_{\mathbb{R}^{d+1}}  \prod^{2\kappa+1}_{\ell=1} \Box^\alpha_{k^{(\ell)}} u_{\ell} \overline{\Box^\alpha_{k} v} dxdt \nonumber\\
& \leq  \|\Box^\alpha_{k} v \|_{L^4_{x,t}}  \|\Box^\alpha_{k^{(1)}} u_1 \|_{L^4_{x,t}} \|\Box^\alpha_{k^{(2)}} u_2 \|_{L^4_{x,t}} \|\Box^\alpha_{k^{(3)}} u_{3} \|_{L^4_{x,t}}  \prod^{2\kappa+1}_{\ell=4} \|\Box^\alpha_{k^{(\ell)}} u_{\ell}\|_{L^\infty_{x,t}} \nonumber\\
& \lesssim   \langle k^{(1)} \rangle^{(d-2)\alpha/ (1-\alpha) } \| \Box^\alpha_{k} v  \|_{V^2_\Delta} \|\Box^\alpha_{k^{(3)}} u  \|_{U^2_\Delta} \|\Box^\alpha_{k^{(1)}} u  \|_{U^2_\Delta} \| \Box^\alpha_{-k^{(2)}} u  \|_{U^2_\Delta}  \nonumber\\
& \quad \times \prod^{2\kappa+1}_{\ell=4} \langle k^{(\ell)} \rangle^{d\alpha/2(1-\alpha)} \|\Box^\alpha_{\pm k^{(\ell)}} u_{\ell}\|_{U^2_\Delta }.
\end{align}
Similar to Case 1.1,
\begin{align}
\mathscr{L}_1 (u,v)
& \lesssim  \sum_{ ((k^{(\ell)})^{2\kappa+1}_{\ell=1} \in \Omega_1, \ k\in \Lambda ((k^{(\ell)})^{2\kappa+1}_{\ell=1}) }   \langle k^{(1)} \rangle^{- 2\alpha /(1-\alpha) + 2\kappa(d\alpha/2-s)/(1-\alpha)}   \langle k  \rangle^{- s/ (1-\alpha)} \| \Box^\alpha_{k} v  \|_{V^2_\Delta}    \nonumber\\
& \quad \times \prod^{2\kappa+1}_{\ell=1} \langle k^{(\ell)} \rangle^{s/ (1-\alpha)} \|\Box^\alpha_{\pm k^{(\ell)}} u \|_{U^2_\Delta } \nonumber\\
& \lesssim   \|   v  \|_{Y^{-s}_\Delta}  \sum_{  k^{(1)},...,k^{(2\kappa+1)}\in \mathbb{Z}^d }
  \prod^{2\kappa+1}_{\ell=1} \langle k^{(\ell)} \rangle^{s/ (1-\alpha)} \|\Box^\alpha_{\pm k^{(\ell)}} u \|_{U^2_\Delta } \nonumber\\
& \lesssim   \|   v  \|_{Y^{-s,\alpha}_\Delta}   \| u \|^{2\kappa+1}_{X^{s,\alpha}_\Delta}.
\end{align}

{\it Case 2. } If $ \langle k_j \rangle  \sim \langle k \rangle  \sim  \langle k^{(2)} \rangle \sim \langle k^{(2)}_j \rangle$, the arguments are similar to Case 1 and we omit the details.

\subsection{Estimates of $\mathscr{L}_0 (u,v)$} \label{estimateL0}

When we estimate $\mathscr{L}_1 (u,v)$, we see that $\# \Lambda ((k^{(\ell)})^{2\kappa+1}_{\ell=1}) \sim 1$, which leads to the summation over all $k\in \Lambda ((k^{(\ell)})^{2\kappa+1}_{\ell=1})$ has no contribution to the regularity index $s>s_\kappa$.
 However, in the case $(k^{(\ell)})^{2\kappa+1}_{\ell=1} \in \Omega_0 \cup \Omega_2$, $\Lambda ((k^{(\ell)})^{2\kappa+1}_{\ell=1})$ is a bit complicated and we have
 \begin{lem} \label{Nonlinearnumberlambda}
We have
\begin{itemize}

\item[\rm (i)] If $K \lesssim \bigvee^{2\kappa+1}_{\ell=1} \langle k^{(\ell)} \rangle^{\alpha/(1-\alpha)}$, then
\begin{align}
\# \Lambda ((k^{(\ell)})^{2\kappa+1}_{\ell=1}) \lesssim   \bigvee^{2\kappa+1}_{\ell=1} \langle k^{(\ell)} \rangle^{d \alpha }.  \label{Numberlambda1}
\end{align}
\item[\rm (ii)] If $K \gg \bigvee^{2\kappa+1}_{\ell=1} \langle k^{(\ell)} \rangle^{\alpha/(1-\alpha)}$, then we have
 \begin{align}
\# \Lambda ((k^{(\ell)})^{2\kappa+1}_{\ell=1}) \lesssim  1 \bigvee  \left( \bigvee^{2\kappa+1}_{\ell=1} \langle k^{(\ell)} \rangle^{d \alpha/(1-\alpha)}/K^{d\alpha} \right).  \label{Numberlambda2}
\end{align}
    \end{itemize}
\end{lem}
Lemma \ref{Nonlinearnumberlambda} is a simple consequence of \eqref{Nonlinearnonzero1} and \eqref{Nonlinearnonzero2} and we omit the details for the proof. In order to estimate $\mathscr{L}_0 (u,v)$, a straightforward idea is to follow the same way as in the estimates of $\mathscr{L}_1 (u,v)$ and use Lemma \ref{Nonlinearnumberlambda} to calculate the number of $\Lambda ((k^{(\ell)})^{2\kappa+1}_{\ell=1})$. Unfortunately, the summation on $k \in \Lambda ((k^{(\ell)})^{2\kappa+1}_{\ell=1})$ will make troubles to the optimal regularity index $s>s_\kappa=d\alpha/2-\alpha/\kappa$. So, we need to use a different way to control $\mathscr{L}_0 (u,v)$ and we use the dyadic decomposition on $\overline{v}$.  We can rewrite $\mathscr{L}_0 (u,v)$ as
\begin{align}
  \mathscr{L}_0 (u,v) = \sum_{(k^{(\ell)})^{2\kappa+1}_{\ell=1} \in \Omega_0, \  j\in \mathbb{Z}_+ } \int_{\mathbb{R}^{d+1}}  \prod^{2\kappa+1}_{\ell=1} \Box^\alpha_{k^{(\ell)}} u_{\ell}  \overline{\triangle_j v} dxdt .
\label{nonlinear01}
\end{align}
So, one needs to estimate the right hand side of  \eqref{nonlinear01} and we prove that
\begin{lem} \label{Nonlinearestim02}

\begin{align}
 |\mathscr{L}_{0} (u,v)| \lesssim   \|   v  \|_{Y^{-s, \alpha}_\Delta}   \| u \|^{2\kappa+1}_{X^{s,\alpha}_\Delta}     \label{nonlinearestim01}
\end{align}
\end{lem}

\begin{lem} \label{Nonlinearestim03}
Let $(k^{(\ell)})^{2\kappa+1}_{\ell=1} \in \Omega_0$. If
$$
 \int_{\mathbb{R}^{d+1}}  \prod^{2\kappa+1}_{\ell=1} \Box^\alpha_{k^{(\ell)}} u_{\ell}  \overline{\triangle_j v} dxdt \neq 0,
$$
then we have
\begin{align} \label{lambda0}
j\in \Lambda((k^{(\ell)})^{2\kappa+1}_{\ell=1} ):=  \left\{j\in \mathbb{Z}_+: \  2^j \lesssim  \bigvee^{2\kappa+1}_{\ell=1} \langle k^{(\ell)} \rangle^{\alpha /(1-\alpha)} \right\}.
\end{align}
\end{lem}
Note here and below in this subsection the notation $\Lambda((k^{(\ell)})^{2\kappa+1}_{\ell=1} )$ is different from \eqref{lambda1} in the previous subsection.   The proof of Lemma \ref{Nonlinearestim03} proceeds in the same way as in Lemma \ref{Nonlinearnonzero} and we omit the details. Hence, we have
\begin{align}
  \mathscr{L}_0 (u,v) = \sum_{  (k^{(\ell)})^{2\kappa+1}_{\ell=1} \in \Omega_0, \ j\in \Lambda((k^{(\ell)})^{2\kappa+1}_{\ell=1} ) } \int_{\mathbb{R}^{d+1}}  \prod^{2\kappa+1}_{\ell=1} \Box^\alpha_{k^{(\ell)}} u_{\ell}  \overline{\triangle_j v} dxdt .
\label{nonlinear02}
\end{align}
One easily sees that
$$
\# \Lambda((k^{(\ell)})^{2\kappa+1}_{\ell=1} ) \lesssim  \max_{1\leq \ell \leq 2\kappa+1}\ln \langle k^{(\ell)} \rangle,
$$
it follows that the summation on $j\in \Lambda((k^{(\ell)})^{2\kappa+1}_{\ell=1} )$ in \eqref{nonlinear02} is bounded by $O(\max_{1\leq \ell \leq 2\kappa+1}\ln \langle k^{(\ell)} \rangle)$, which is much less than that of any $\alpha$-decomposition.  This is the main reason why we try to use the dyadic decomposition with respect to $\overline{v}$.

We can assume that $\max_{1\leq \ell \leq 2\kappa+1} |k^{(\ell)}| \gg 1$. In the opposite case all $u_\ell$ and $\overline{v}$ have lower frequency,  \eqref{nonlinearestim01} can be easily obtained.  As in \eqref{Numberorder}, we can assume that
\begin{align}
  |k^{(1)}| \geq |k^{(3)}|\geq ...\geq |k^{(2\kappa+1)}|, \ \ |k^{(2)}| \geq |k^{(4)}|\geq ...\geq |k^{(2\kappa)}|. \label{Numberorder0}
\end{align}

We divide the proof of Lemma \ref{Nonlinearestim02} in to the following three cases.

{\it Case 1.} We consider $| k^{(1)}| =\max _{1\leq \ell \leq 2\kappa+1} | k^{(\ell)}| \gg 1$.
By \eqref{Numberorder0} and $(k^{(1)},...,k^{(2\kappa+1)}) \in \Omega_0 $,  we have
\begin{align}
   \langle k^{(1)} \rangle \sim  \langle k^{(2)}  \rangle, \ or  \    \langle k^{(1)} \rangle \sim \langle k^{(3)}  \rangle.
\end{align}
If not, then $|k^{(\ell)}| \ll |k^{(1)}|$ for all $\ell=2,...,2\kappa+1$. Assume that $|k^{(1)}_1| \geq |k^{(1)}_i|$ for $i=2,...,d$. Then
$$
2^j \geq  \langle k^{(1)} \rangle^{\alpha/(1-\alpha)} (| k^{(1)}_1| -C) -  \sum^{2\kappa+1}_{\ell=2} \langle k^{(\ell)} \rangle^{\alpha/(1-\alpha)} (| k^{(\ell)}_1| +C)  \gtrsim   \langle k^{(1)} \rangle^{1/(1-\alpha)} \gg \bigvee^{2\kappa+1}_{\ell=1} \langle k^{(\ell)} \rangle^{\alpha /(1-\alpha)}.
$$
A contradiction.  The case $\langle k^{(1)} \rangle \sim \langle k^{(3)}  \rangle$ is similar to the case $\langle k^{(1)} \rangle \sim \langle k^{(2)}  \rangle$, it suffices to consider the case $\langle k^{(1)} \rangle \sim \langle k^{(2)}  \rangle$.

{\it Case 1.1.}  We consider the case $ \langle k^{(1)} \rangle \sim  \langle k^{(2)} \rangle  \gg  \langle k^{(\ell)} \rangle$, $\ell=3,...,2\kappa+1$.
By H\"older's inequality, bilinear estimates, we have
\begin{align}
& \int_{\mathbb{R}^{d+1}}  \prod^{2\kappa+1}_{\ell=1} \Box^\alpha_{k^{(\ell)}} u_{\ell} \overline{\triangle_j v} dxdt \nonumber\\
& \leq \|\overline{\triangle_j v} \Box^\alpha_{k^{(1)}} u_{1} \|_{L^2_{x,t}} \|\Box^\alpha_{k^{(2)}} u_2  \Box^\alpha_{k^{(3)}} u_{3} \|_{L^2_{x,t}}  \prod^{2\kappa+1}_{\ell=4} \|\Box^\alpha_{k^{(\ell)}} u_{\ell}\|_{L^\infty_{x,t}} \nonumber\\
& \lesssim \ln \langle k^{(1)}\rangle \  \langle k^{(1)} \rangle^{-1/(1-\alpha)}  2^{j(d-1)/2} \langle k^{(3)} \rangle^{(d-1)\alpha/2(1-\alpha)} \| \triangle_j v \|_{V^2_\Delta} \|\Box^\alpha_{k^{(1)}} u  \|_{U^2_\Delta} \|\Box^\alpha_{-k^{(2)}} u  \|_{U^2_\Delta} \| \Box^\alpha_{k^{(3)}} u  \|_{U^2_\Delta}  \nonumber\\
& \quad \times \prod^{2\kappa+1}_{\ell=4} \langle k^{(\ell)} \rangle^{d\alpha/2(1-\alpha)} \|\Box^\alpha_{\pm k^{(\ell)}} u_{\ell}\|_{U^2_\Delta } \nonumber\\
& \lesssim  \ln \langle k^{(1)}\rangle \  \langle k^{(1)} \rangle^{-(1+2s)/(1-\alpha)}  2^{j((d-1)/2 + s + d(1-\alpha)/2)} 2^{j(-s-d(1-\alpha)/2)} \|\triangle_j v  \|_{V^2_\Delta}   \nonumber\\
& \quad \times  \langle k^{(3)} \rangle^{(\kappa(d\alpha/2 -s) -\alpha/2)/(1-\alpha)}\langle k^{(4)} \rangle ^{(\kappa-1)(d\alpha/2-s)/(1-\alpha)} \prod^{2\kappa+1}_{\ell=1} \langle k^{(\ell)} \rangle^{s/(1-\alpha)} \|\Box^\alpha_{\pm k^{(\ell)}} u_{\ell}\|_{U^2_\Delta }. \label{L01}
\end{align}
Since
$ \kappa (d\alpha/2-s)  - \alpha/2  > 0$
and $|k^{(3)}|, |k^{(4)}| \leq |k^{(1)}|$ in the summation of $\mathscr{L}_0 (u,v)$, noticing that
\begin{align}
\sum_{ j\in \Lambda ((k^{(\ell)})^{2\kappa+1}_{\ell=1}) }  2^{j((d-1)/2 + s + d(1-\alpha)/2)} \lesssim   \langle k^{(1)} \rangle^{\alpha ((d-1)/2 + s + d(1-\alpha)/2))/(1-\alpha)}, \label{dyadicsum}
\end{align}
 we have
\begin{align}
\mathscr{L}_0 (u,v)
& \lesssim  \sum_{ ((k^{(\ell)})^{2\kappa+1}_{\ell=1} \in \Omega_0, \ j\in \Lambda ((k^{(\ell)})^{2\kappa+1}_{\ell=1}) } \ln \langle k^{(1)}\rangle \  \langle k^{(1)} \rangle^{(- 1-2s-\alpha/2  + (2\kappa-1)(d\alpha/2-s))/(1-\alpha)} \nonumber\\
& \quad \times \prod^{2\kappa+1}_{\ell=1} \langle k^{(\ell)} \rangle^{s/ (1-\alpha)} \|\Box^\alpha_{\pm k^{(\ell)}} u \|_{U^2_\Delta } \ 2^{j((d-1)/2 + s + d(1-\alpha)/2)} \|v\|_{Y^{-s-d(1-\alpha)/2, 1}_\Delta}  \nonumber\\
& \lesssim  \|   v  \|_{Y^{-s-d(1-\alpha)/2, 1}_\Delta}  \sum_{  (k^{(\ell)}) \in \mathbb{Z}^{d(2\kappa+1)} } \ln \langle k^{(1)}\rangle \  \langle k^{(1)} \rangle^{((2\kappa+1-\alpha)(d\alpha/2-s) -1-\alpha)/(1-\alpha) } \nonumber\\
& \quad\quad  \times \prod^{2\kappa+1}_{\ell=1} \langle k^{(\ell)} \rangle^{s/ (1-\alpha)} \|\Box^\alpha_{\pm k^{(\ell)}} u \|_{U^2_\Delta } \nonumber\\
& \lesssim   \|   v  \|_{Y^{-s-d(1-\alpha)/2, 1}_\Delta}   \| u \|^{2\kappa+1}_{X^{s,\alpha}_\Delta}.
\end{align}

{\it Case 1.2.} We consider the case $ \langle k^{(1)} \rangle \sim  \langle k^{(2)} \rangle  \sim  \langle k^{(3)} \rangle$, $\kappa\geq 2$.
By H\"older's inequality, bilinear and the Strichartz estimates we have
\begin{align}
& \int_{\mathbb{R}^{d+1}}  \prod^{2\kappa+1}_{\ell=1} \Box^\alpha_{k^{(\ell)}} u_{\ell} \overline{\triangle_j v} dxdt \nonumber\\
& \leq \|\overline{\triangle_j v} \Box^\alpha_{k^{(1)}} u_{1} \|_{L^2_{x,t}} \|\Box^\alpha_{k^{(2)}} u_2\|_{L^4_{x,t}}  \|\Box^\alpha_{k^{(3)}} u_{3} \|_{L^4_{x,t}}  \prod^{2\kappa+1}_{\ell=4} \|\Box^\alpha_{k^{(\ell)}} u_{\ell}\|_{L^\infty_{x,t}} \nonumber\\
& \lesssim \ln \langle k^{(1)}\rangle \  \langle k^{(1)} \rangle^{-1/2(1-\alpha)}  2^{j(d-1)/2} \langle k^{(2)} \rangle^{(d\alpha/4-\alpha/2)/ (1-\alpha)} \langle k^{(3)} \rangle^{(d\alpha/4-\alpha/2)/ (1-\alpha)}\| \triangle_j v \|_{V^2_\Delta} \|\Box^\alpha_{k^{(1)}} u  \|_{U^2_\Delta} \nonumber\\
& \quad \times \|\Box^\alpha_{-k^{(2)}} u  \|_{U^2_\Delta} \| \Box^\alpha_{k^{(3)}} u  \|_{U^2_\Delta}  \prod^{2\kappa+1}_{\ell=4} \langle k^{(\ell)} \rangle^{d\alpha/2(1-\alpha)} \|\Box^\alpha_{\pm k^{(\ell)}} u_{\ell}\|_{U^2_\Delta } \nonumber\\
& \lesssim  2^{j((d-1)/2 + s + d(1-\alpha)/2)} 2^{j(-s-d(1-\alpha)/2)} \|\triangle_j v  \|_{V^2_\Delta}   \nonumber\\
& \quad \times \ln \langle k^{(1)}\rangle \  \langle k^{(1)} \rangle^{(-1/2-2s-\alpha + (2\kappa-1)(d\alpha/2-s))/(1-\alpha)}  \prod^{2\kappa+1}_{\ell=1} \langle k^{(\ell)} \rangle^{s/(1-\alpha)} \|\Box^\alpha_{\pm k^{(\ell)}} u_{\ell}\|_{U^2_\Delta }.
\end{align}
In the case $\kappa\geq 2$, noticing \eqref{dyadicsum}, we have
 we have
\begin{align}
\mathscr{L}_0 (u,v)
& \lesssim  \sum_{ ((k^{(\ell)})^{2\kappa+1}_{\ell=1} \in \Omega_0, \ j\in \Lambda ((k^{(\ell)})^{2\kappa+1}_{\ell=1}) }  \ln \langle k^{(1)}\rangle \  \langle k^{(1)} \rangle^{(- 1-2s-\alpha  + (2\kappa-1)(d\alpha/2-s))/(1-\alpha)} \nonumber\\
& \quad \times \prod^{2\kappa+1}_{\ell=1} \langle k^{(\ell)} \rangle^{s/ (1-\alpha)} \|\Box^\alpha_{\pm k^{(\ell)}} u \|_{U^2_\Delta } 2^{j((d-1)/2 + s + d(1-\alpha)/2)} \|v\|_{Y^{-s-d(1-\alpha)/2, 1}_\Delta}  \nonumber\\
& \lesssim  \|   v  \|_{Y^{-s-d(1-\alpha)/2, 1}_\Delta}  \sum_{  (k^{(\ell)}) \in \mathbb{Z}^{d(2\kappa+1)} } \ln \langle k^{(1)}\rangle \  \langle k^{(1)} \rangle^{((2\kappa+1-\alpha)(d\alpha/2-s) -1/2-3\alpha/2)/(1-\alpha) } \nonumber\\
& \quad \quad \times  \prod^{2\kappa+1}_{\ell=1} \langle k^{(\ell)} \rangle^{s/ (1-\alpha)} \|\Box^\alpha_{\pm k^{(\ell)}} u \|_{U^2_\Delta } \nonumber\\
& \lesssim   \|   v  \|_{Y^{-s-d(1-\alpha)/2, 1}_\Delta}   \| u \|^{2\kappa+1}_{X^{s,\alpha}_\Delta}.
\end{align}
If $\kappa=1$, we will give the proof in the next Section.

{\it Case 1.3.} If $ \langle k^{(1)} \rangle \sim  \langle k^{(2)} \rangle  \sim  \langle k^{(4)} \rangle$, $\kappa\geq 2$,  one can use the same way as Case 1.2 to show the result and we omit the details.

{\it Case 2.} We consider $| k^{(2)}| =\max _{1\leq \ell \leq 2\kappa+1} | k^{(\ell)}| \gg 1$. The argument is similar to Case 1.

\subsection{Estimates of $\mathscr{L}_2 (u,v)$} \label{EstimateL2}

If $(k^{(\ell)})^{2\kappa+1}_{\ell=1} \in \Omega_2$, we see that there exists $0<\theta<1$ verifying ${K}=(\bigvee^{2\kappa+1}_{\ell=1}\langle k^{(\ell)} \rangle^{1/(1-\alpha)}) ^{1-\theta+\alpha\theta}$. Put
\begin{equation}\label{alpha-tilde}
  \widetilde{\alpha}=\frac{\alpha}{1-\theta+\alpha\theta}.
\end{equation}
Using the $\widetilde{\alpha}$-decomposition to $\overline{v}$,   one can rewrite $\mathscr{L}_2 (u,v)$ as

\begin{align}
  \mathscr{L}_2 (u,v) = \sum_{ k\in \mathbb{Z}^d,  \ (k^{(\ell)}) \in \Omega_2 } \int_{\mathbb{R}^{d+1}}  \prod^{2\kappa+1}_{\ell=1} \Box^\alpha_{k^{(\ell)}} u_{\ell}  \overline{\Box^{\widetilde{\alpha}}_{k} v} dxdt .
\label{nonlinear231}
\end{align}
So, one needs to control the right hand side of  \eqref{nonlinear231} and we prove that
\begin{lem} \label{Nonlinearestim1}

\begin{align}
 |\mathscr{L}_{2} (u,v)| \lesssim   \|   v  \|_{Y^{-s, \alpha}_\Delta}   \| u \|^{2\kappa+1}_{X^{s,\alpha}_\Delta}     \label{nonlinearestim2-1}
\end{align}
\end{lem}

First, we have

\begin{lem} \label{Nonlinearnonzero-2}
If $\mathscr{L}_{2} (u,v)$ in \eqref{nonlinear231} does not equal to zero, we must have
\begin{align}
 \langle k \rangle^{\widetilde{\alpha}/(1-\widetilde{\alpha})} (k_j -C) \leq  \sum^{2\kappa+1}_{\ell =1} \langle k^{(\ell)} \rangle^{\alpha/(1-\alpha)} (k^{(\ell)}_j + C),   \label{Nonlinearnonzero2-1}\\
 \langle k \rangle^{\widetilde{\alpha}/(1-\widetilde{\alpha})} (k_j + C) \geq  \sum^{2\kappa+1}_{\ell =1} \langle k^{(\ell)} \rangle^{\alpha/(1-\alpha)} (k^{(\ell)}_j - C).   \label{Nonlinearnonzero2-2}
\end{align}
\end{lem}
{\bf Proof.} See Lemma \ref{Nonlinearnonzero}. $\hfill\Box$

The idea to use the $\widetilde{\alpha}$-decomposition is similar to the cases as in handling $\mathscr{L}_{0} (u,v)$,  one needs to remove the summation on $k\in \mathbb{Z}^d$ in the right hand side of \eqref{nonlinear231}. For convenience, we write
\begin{align}
\Lambda((k^{(\ell)})^{2\kappa+1}_{\ell=1}) = \{k \in \mathbb{Z}^{d} : \eqref{Nonlinearnonzero2-1} \ {\rm and}  \ \eqref{Nonlinearnonzero2-2} \ {\rm are \ satisfied}\}. \label{lambda2-1}
\end{align}
Note here $\Lambda((k^{(\ell)})^{2\kappa+1}_{\ell=1})$ has different meaning as in the previous two subsections. So, we can further rewrite $\mathscr{L}_2 (u,v)$ as
\begin{align}
  \mathscr{L}_1 (u,v) = \sum_{ ((k^{(\ell)})^{2\kappa+1}_{\ell=1} \in \Omega_2, \ k\in \Lambda ((k^{(\ell)})^{2\kappa+1}_{\ell=1}) } \int_{\mathbb{R}^{d+1}}  \prod^{2\kappa+1}_{\ell=1} \Box^\alpha_{k^{(\ell)}} u_{\ell}  \overline{\Box^{\widetilde{\alpha}}_{k} v} dxdt .
\label{nonlinear2-11}
\end{align}

\begin{lem} \label{Nonlinearnonzero-2-3}
Let $ ((k^{(\ell)})^{2\kappa+1}_{\ell=1} \in \Omega_2$.  We have $\# \Lambda((k^{(\ell)})^{2\kappa+1}_{\ell=1}) \lesssim 1$  and $K\sim \langle k \rangle^{1/(1-\widetilde{\alpha})} $ for any $k\in \Lambda((k^{(\ell)})^{2\kappa+1}_{\ell=1})$.

\end{lem}
{\bf Proof.} First, we show that $K\sim \langle k \rangle^{1/(1-\widetilde{\alpha})} $ for any $k\in \Lambda((k^{(\ell)})^{2\kappa+1}_{\ell=1})$.
We can assume that $K_1= \max_{1\leq i\leq 2\kappa+1} K_i =K$.  Since $((k^{(\ell)})^{2\kappa+1}_{\ell=1} \in \Omega_2$, we see that
$$
\left|\sum^{2\kappa+1}_{\ell =1} \langle k^{(\ell)} \rangle^{\alpha/(1-\alpha)} (k^{(\ell)}_1 \pm C)\right| \sim K.
$$
Using \eqref{Nonlinearnonzero2-1} and  \eqref{Nonlinearnonzero2-2}, we get $K\sim \langle k \rangle^{1/(1-\widetilde{\alpha})} $.  Let $k=(k_1, \overline{k}), \ k'=(k'_1, \overline{k}) \in \Lambda((k^{(\ell)})^{2\kappa+1}_{\ell=1})$. In view of mean value Theorem together with \eqref{Nonlinearnonzero2-1} and  \eqref{Nonlinearnonzero2-2}, one has that
$$
|k_1-k'_1|  \lesssim  K^{-\widetilde{\alpha}}  \left(K^{\widetilde{\alpha}} + \bigvee^{2\kappa+1}_{\ell=1} \langle k^{(\ell)} \rangle^{\alpha/(1-\alpha)} \right)  \lesssim 1.
$$

We consider the estimates of $\mathscr{L}_2 (u,v)$.  As before,  we can assume that
\begin{align}
  |k^{(1)}| \geq |k^{(3)}|\geq ...\geq |k^{(2\kappa+1)}|, \ \ |k^{(2)}| \geq |k^{(4)}|\geq ...\geq |k^{(2\kappa)}|. \label{Numberorder2-0}
\end{align}
For $(k^{(1)},...,k^{(2\kappa+1)}) \in \Omega_2 $, we see that $ |k^{(1)}|\vee  |k^{(2)}| =\max_{1\leq \ell \leq 2\kappa+1} |k^{(\ell)}|$.

{\it Case 1.}  We assume that $ |k^{(1)}| =\max_{1\leq \ell \leq 2\kappa+1} |k^{(\ell)}|$.  It is easy to see that
\begin{align}
  \langle k^{(1)} \rangle \sim  \langle k^{(2)}  \rangle, \ or  \    \langle k^{(1)} \rangle \sim \langle k^{(3)}  \rangle.
\end{align}
We can further assume that $\langle k^{(1)} \rangle \sim  \langle k^{(2)}  \rangle$ and the case $\langle k^{(1)} \rangle \sim \langle k^{(3)}  \rangle$ can be handled by using a similar way.
 There exists $j\in \{1,...,d\}$ such that $\langle k^{(1)} \rangle \sim  \langle k^{(1)}_j \rangle.$
We can assume that $j=1$.

{\it Case 1.1.} We consider the case $ |k^{(\ell)}| \ll |k^{(1)}|$, $\ell =3,...,2\kappa+1$. We have $\langle k^{(1)} \rangle^{1/(1-\alpha)} \gg K \sim \langle k \rangle^{1/(1-\widetilde{\alpha})}$.  In $\mathscr{L}_2(u,v)$, using H\"older's, bilinear and Strichartz' estimates, we can bound
\begin{align}
& \int_{\mathbb{R}^{d+1}}  \prod^{2\kappa+1}_{\ell=1} \Box^\alpha_{k^{(\ell)}} u_{\ell} \overline{\Box^{\widetilde{\alpha}}_{k} v} dxdt \nonumber\\
& \leq \|\overline{\Box^{\widetilde{\alpha}}_{k} v} \Box^\alpha_{k^{(1)}} u_{1} \|_{L^2_{x,t}} \|\Box^\alpha_{k^{(2)}} u_2  \Box^\alpha_{k^{(3)}} u_{3} \|_{L^2_{x,t}}  \prod^{2\kappa+1}_{\ell=4} \|\Box^\alpha_{k^{(\ell)}} u_{\ell}\|_{L^\infty_{x,t}} \nonumber\\
& \lesssim  \ln \langle k^{(1)}\rangle \  \langle k^{(1)} \rangle^{-1/(1-\alpha)}  \langle k  \rangle^{(d-1)\widetilde{\alpha}/2(1-\widetilde{\alpha})} \langle k^{(3)} \rangle^{(d-1)\alpha/2(1-\alpha)} \| \Box^{\widetilde{\alpha}}_{k} v  \|_{V^2_\Delta} \|\Box^\alpha_{k^{(1)}} u  \|_{U^2_\Delta}  \nonumber\\
& \quad \times \| \Box^\alpha_{-k^{(2)}} u  \|_{U^2_\Delta} \|\Box^\alpha_{k^{(3)}} u  \|_{U^2_\Delta} \prod^{2\kappa+1}_{\ell=4} \langle k^{(\ell)} \rangle^{d\alpha/2(1-\alpha)} \|\Box^\alpha_{\pm k^{(\ell)}} u_{\ell}\|_{U^2_\Delta } \nonumber\\
& \lesssim \ln \langle k^{(1)}\rangle \  \langle k^{(1)} \rangle^{(-1-2s)/(1-\alpha)}  \langle k  \rangle^{((d-1)\widetilde{\alpha}/2 +s+ d(\widetilde{\alpha}-\alpha)/2)/(1-\widetilde{\alpha})}   \langle  k  \rangle^{(-s-d(\widetilde{\alpha}-\alpha)/2) /(1-\widetilde{\alpha})}\| \Box^\alpha_{k} v  \|_{V^2_\Delta}   \nonumber\\
& \quad \times   \langle k^{(3)} \rangle^{( \kappa (d\alpha/2-s)-\alpha/2)/(1-\alpha)} \langle k^{(4)} \rangle^{( \kappa-1) (d\alpha/2-s)/(1-\alpha)}\prod^{2\kappa+1}_{\ell=1} \langle k^{(\ell)} \rangle^{s/(1-\alpha)} \|\Box^\alpha_{\pm k^{(\ell)}} u_{\ell}\|_{U^2_\Delta }. \label{nonlinear2-1}
\end{align}
 Notice that $\langle k\rangle^{1/(1-\widetilde{\alpha})} \lesssim \langle k^{(1)}\rangle^{(1-\theta+\alpha\theta)/(1-{\alpha})} $ and
\begin{align}
& \kappa (d\alpha/2-s) -\alpha/2 \geq 0, \label{nonlinear2-2}\\
&  \langle k  \rangle^{((d-1)\widetilde{\alpha}/2 +s+ d(\widetilde{\alpha}-\alpha)/2)/(1-\widetilde{\alpha})} \lesssim \langle k^{(1)}  \rangle^{((d-1)\alpha/2 +s + \theta(1-\alpha)(d\alpha/2-s))/(1-\alpha)}, \label{nonlinear2-3}\\
& (2\kappa+\theta(1-\alpha))(d\alpha/2-s) <  1+ \alpha. \label{nonlinear2-4}
\end{align}
By Lemma \ref{Nonlinearnonzero-2-3},
we have $\# \Lambda(k^{(1)},...,k^{(2\kappa+1)}) \lesssim 1$ if $(k^{(1)},...,k^{(2\kappa+1)}) \in \Omega_2$. We have from \eqref{nonlinear2-1}--\eqref{nonlinear2-4} that
\begin{align}
\mathscr{L}_2 (u,v)
& \lesssim  \sum_{ ((k^{(\ell)})^{2\kappa+1}_{\ell=1} \in \Omega_2, \ k\in \Lambda ((k^{(\ell)})^{2\kappa+1}_{\ell=1}) } \ln \langle k^{(1)}\rangle \  \langle k^{(1)} \rangle^{(-1-\alpha+ (2\kappa+\theta(1-\alpha))(d\alpha/2-s))/(1-\alpha)}    \nonumber\\
& \quad \times \langle  k  \rangle^{(-s-d(\widetilde{\alpha}-\alpha)/2) /(1-\widetilde{\alpha})} \| \Box^\alpha_{k} v  \|_{V^2_\Delta} \prod^{2\kappa+1}_{\ell=1} \langle k^{(\ell)} \rangle^{s/ (1-\alpha)} \|\Box^\alpha_{\pm k^{(\ell)}} u \|_{U^2_\Delta } \nonumber\\
& \lesssim   \|   v  \|_{Y^{-s-d(\widetilde{\alpha}-\alpha)/2,\widetilde{\alpha}}_\Delta}   \| u \|^{2\kappa+1}_{X^{s,\alpha}_\Delta}. \label{L2est}
\end{align}

{\it Case 1.2.} We consider the case $|k^{(1)}| \sim |k^{(2)}| \sim |k^{(3)}| $ and $\kappa\geq 2$; or $|k^{(1)}| \sim |k^{(2)}| \sim |k^{(4)}| $ and $\kappa\geq 2$. It suffices to consider the case $|k^{(1)}| \sim |k^{(2)}| \sim |k^{(3)}| $ and $\kappa\geq 2$. We have $\langle k^{(1)} \rangle^{1/\alpha} \gg K$.  In $\mathscr{L}_2(u,v)$, applying H\"older's, bilinear and Strichartz' estimates, we have
\begin{align}
& \int_{\mathbb{R}^{d+1}}  \prod^{2\kappa+1}_{\ell=1} \Box^\alpha_{k^{(\ell)}} u_{\ell} \overline{\Box^{\widetilde{\alpha}}_{k} v} dxdt \nonumber\\
& \leq \|\overline{\Box^{\widetilde{\alpha}}_{k} v} \Box^\alpha_{k^{(1)}} u_{1} \|_{L^2_{x,t}} \|\Box^\alpha_{k^{(2)}} u_2  \Box^\alpha_{k^{(3)}} u_{3} \|_{L^2_{x,t}}  \prod^{2\kappa+1}_{\ell=4} \|\Box^\alpha_{k^{(\ell)}} u_{\ell}\|_{L^\infty_{x,t}} \nonumber\\
& \lesssim \ln \langle k^{(1)}\rangle \  \langle k^{(1)} \rangle^{-1/2(1-\alpha)}  \langle k \rangle^{(d-1)\widetilde{\alpha}/2(1-\widetilde{\alpha})} \langle k^{(2)} \rangle^{(d\alpha/4-\alpha/2)/2(1-\alpha)} \langle k^{(3)} \rangle^{(d\alpha/4-\alpha/2)/2(1-\alpha)}  \nonumber\\
& \quad \times \| \Box^{\widetilde{\alpha}}_{k} v  \|_{V^2_\Delta} \prod^3_{\ell=1} \|\Box^\alpha_{\pm k^{(\ell)}} u  \|_{U^2_\Delta} \prod^{2\kappa+1}_{\ell=4} \langle k^{(\ell)} \rangle^{d\alpha/2(1-\alpha)} \|\Box^\alpha_{\pm k^{(\ell)}} u\|_{U^2_\Delta } \nonumber\\
& \lesssim  \ln \langle k^{(1)}\rangle \   \langle k  \rangle^{((d-1)\widetilde{\alpha}/2 +s+ d(\widetilde{\alpha}-\alpha)/2)/(1-\widetilde{\alpha})}   \langle  k  \rangle^{(-s-d(\widetilde{\alpha}-\alpha)/2) /(1-\widetilde{\alpha})}\| \Box^\alpha_{k} v  \|_{V^2_\Delta}   \nonumber\\
& \quad \times   \langle k^{(1)} \rangle^{(( 2\kappa-1) (d\alpha/2-s)-1/2-\alpha-2s)/(1-\alpha)}\prod^{2\kappa+1}_{\ell=1} \langle k^{(\ell)} \rangle^{s/(1-\alpha)} \|\Box^\alpha_{\pm k^{(\ell)}} u\|_{U^2_\Delta }. \label{nonlinear2-5}
\end{align}
Since $\kappa\geq 2$, we have from \eqref{nonlinear2-3} and \eqref{nonlinear2-5} that
\begin{align}
\mathscr{L}_2 (u,v)
& \lesssim  \sum_{ ((k^{(\ell)})^{2\kappa+1}_{\ell=1} \in \Omega_2, \ k\in \Lambda ((k^{(\ell)})^{2\kappa+1}_{\ell=1}) }   \ln \langle k^{(1)}\rangle \   \langle k^{(1)} \rangle^{(-1/2-3\alpha/2 + (2\kappa+\theta(1-\alpha))(d\alpha/2-s))/(1-\alpha)}    \nonumber\\
& \quad \times \langle  k  \rangle^{(-s-d(\widetilde{\alpha}-\alpha)/2) /(1-\widetilde{\alpha})} \| \Box^\alpha_{k} v  \|_{V^2_\Delta} \prod^{2\kappa+1}_{\ell=1} \langle k^{(\ell)} \rangle^{s/ (1-\alpha)} \|\Box^\alpha_{\pm k^{(\ell)}} u \|_{U^2_\Delta } \nonumber\\
& \lesssim   \|   v  \|_{Y^{-s-d(\widetilde{\alpha}-\alpha)/2,\widetilde{\alpha}}_\Delta}   \| u \|^{2\kappa+1}_{X^{s,\alpha}_\Delta}.
\end{align}
We will consider the case $\kappa=1$ in the next subsection.

\subsection{Estimates of $\mathscr{L}_0(u,v)$ and $\mathscr{L}_2(u,v)$ for $\kappa=1$}

 {\bf Step 1.}  (Estimates of $\mathscr{L}_0(u,v)$)  In order to finish the estimates of $\mathscr{L}_0(u,v)$, we need to consider the following case

$$(k^{(1)}, k^{(2)}, k^{(3)}) \in \Omega_0,  \  \  \langle k^{(1)} \rangle \sim  \langle k^{(2)} \rangle  \sim  \langle k^{(3)} \rangle, \ \ \kappa=1.$$

It suffices to consider the case $|k^{(1)}|\gg 1.$ Now we choose $|k^{(1)}_j| = \max_{1\leq i\leq d} |k^{(1)}_i|$, say $|k^{(1)}_1| = \max_{1\leq i\leq d} |k^{(1)}_i|$.  Since
\begin{align} \label{K1}
K_1 = \left| \sum^3_{\ell =1}\langle k^{(\ell)} \rangle^{\alpha/(1-\alpha)} k^{(\ell)}_1  \right|  \lesssim \bigvee^3_{\ell=1} \langle k^{(\ell)} \rangle^{\alpha/(1-\alpha)},
\end{align}
we must have $|k^{(1)}_1| \sim |k^{(2)}_1|$ or $|k^{(1)}_1| \sim |k^{(3)}_1|$. If not, we have $|k^{(2)}_1| \vee |k^{(3)}_1| \ll |k^{(1)}_1|$, then
$$
K_1 \geq \langle k^{(1)} \rangle^{\alpha/(1-\alpha)} |k^{(1)}_1| - \left| \sum^3_{\ell =2}\langle k^{(\ell)} \rangle^{\alpha/(1-\alpha)} k^{(\ell)}_1  \right|  \gtrsim  \langle k^{(1)} \rangle^{1/(1-\alpha)} \gg \bigvee^3_{\ell=1} \langle k^{(\ell)} \rangle^{\alpha/(1-\alpha)}.
$$
A contradiction.

{\it Case 1.} We consider the case $|k^{(1)}_1| = \max_{1\leq i\leq d} |k^{(1)}_i| \sim |k^{(2)}_1|$.

{\it Case 1.1.}  $ |k^{(1)}_1| = \max_{1\leq i\leq d} |k^{(1)}_i| \sim |k^{(2)}_1| \gg |k^{(3)}_1|$.  One can proceeds in the same way as in
Case 1.1 in Section \ref{estimateL0} to obtain the result, the details are omitted.

{\it Case 1.2.}  $ |k^{(1)}_1| = \max_{1\leq i\leq d} |k^{(1)}_i| \sim |k^{(2)}_1| \sim |k^{(3)}_1|$. First, we conclude that $k^{(1)}_1, \ k^{(2)}_1, \ k^{(3)}_1$ cannot have the same signs, i.e., we must have one of the following cases:
\begin{align}
& k^{(1)}_1>0, \ \ k^{(2)}_1<0, \ k^{(3)}_1<0; \ \ or \ k^{(1)}_1<0, \ \ k^{(2)}_1>0, \ k^{(3)}_1>0, \label{+--}\\
& k^{(3)}_1>0, \ \ k^{(1)}_1<0, \ k^{(2)}_1<0; \ \ or \ k^{(3)}_1<0, \ \ k^{(1)}_1>0, \ k^{(2)}_1>0, \label{--+}\\
& k^{(2)}_1>0, \ \ k^{(1)}_1<0, \ k^{(3)}_1<0; \ \ or \ k^{(2)}_1<0, \ \ k^{(1)}_1>0, \ k^{(3)}_1>0. \label{-+-}
\end{align}
If \eqref{+--} holds, by H\"older's inequality,  we have
\begin{align} \label{+--1}
  \int_{\mathbb{R}^{d+1}}  \prod^{3}_{\ell=1} \Box^\alpha_{k^{(\ell)}} u_{\ell} \overline{\triangle_j v} dxdt
& \leq \|\overline{\triangle_j v} \Box^\alpha_{k^{(2)}} \overline{u}  \|_{L^2_{x,t}} \|\Box^\alpha_{k^{(1)}} u \Box^\alpha_{k^{(3)}} u \|_{L^2_{x,t}}
\end{align}
To control the right hand side of \eqref{+--1}, notocing that
\begin{align}
|\langle k^{(1)} \rangle^{\alpha/(1-\alpha)} k^{(1)}_1 -  \langle k^{(\ell)} \rangle^{\alpha/(1-\alpha)} k^{(3)}_1| \gtrsim  \langle k^{(1)} \rangle^{1/(1-\alpha)}  \label{+--2}
\end{align}
and applying the bilinear estimates, we have

\begin{align}
& \int_{\mathbb{R}^{d+1}}  \prod^{3}_{\ell=1} \Box^\alpha_{k^{(\ell)}} u_{\ell} \overline{\triangle_j v} dxdt \nonumber\\
& \lesssim \ln \langle k^{(1)}\rangle \ \langle k^{(1)} \rangle^{-1/(1-\alpha)}  2^{j(d-1)/2} \langle k^{(3)} \rangle^{(d-1)\alpha/2(1-\alpha)} \| \triangle_j v \|_{V^2_\Delta} \|\Box^\alpha_{k^{(1)}} u  \|_{U^2_\Delta} \|\Box^\alpha_{-k^{(2)}} u  \|_{U^2_\Delta} \| \Box^\alpha_{k^{(3)}} u  \|_{U^2_\Delta}  \nonumber\\
& \lesssim \ln \langle k^{(1)}\rangle \  \langle k^{(1)} \rangle^{-(1+2s)/(1-\alpha)}  2^{j((d-1)/2 + s + d(1-\alpha)/2)} 2^{j(-s-d(1-\alpha)/2)} \|\triangle_j v  \|_{V^2_\Delta}   \nonumber\\
& \quad \times  \langle k^{(3)} \rangle^{( (d\alpha/2 -s) -\alpha/2)/(1-\alpha)}  \prod^{3}_{\ell=1} \langle k^{(\ell)} \rangle^{s/(1-\alpha)} \|\Box^\alpha_{\pm k^{(\ell)}} u_{\ell}\|_{U^2_\Delta }. \label{kappa=1-1}
\end{align}
This reduces to the estimate \eqref{L01}.

If \eqref{+--} holds, $k^{(1)}$ and $k^{(3)}$ have equal positions, this case is the same as the above.

If \eqref{-+-} holds,  we have
\begin{align} \label{-+-1}
  \int_{\mathbb{R}^{d+1}}  \prod^{3}_{\ell=1} \Box^\alpha_{k^{(\ell)}} u_{\ell} \overline{\triangle_j v} dxdt
& \leq \|\overline{\triangle_j v} \Box^\alpha_{k^{(3)}} \overline{u}  \|_{L^2_{x,t}} \|\Box^\alpha_{k^{(1)}} u \Box^\alpha_{k^{(2)}} \overline{u} \|_{L^2_{x,t}}
\end{align}
By
$$|\langle k^{(1)} \rangle^{\alpha/(1-\alpha)} k^{(1)}_1 +  \langle k^{(2)} \rangle^{\alpha/(1-\alpha)} k^{(2)}_1| \gtrsim  \langle k^{(3)} \rangle^{1/(1-\alpha)} -K_1 \gtrsim  \langle k^{(1)} \rangle^{1/(1-\alpha)}     $$
and the bilinear estimates in Corollary \ref{bilinear-alpha-schrodinger-222}, we still have \eqref{kappa=1-1}.  Repeating the procedures as above, we have the result, as desired.

{\bf Step 2.} (Estimates of $\mathscr{L}_2(u,v)$)  Now let us connect the proof with Section \ref{EstimateL2} and we use the same notations as in Section \ref{EstimateL2}.  We need to consider the following case
 $$
 (k^{(1)}, k^{(2)}, k^{(3)}) \in \Omega_2, \ \ |k^{(1)}| \sim |k^{(2)}| \sim |k^{(3)}| \gg 1, \ \ \kappa=1.
 $$
Assume that $|k^{(1)}_1| =\max_{1\leq j\leq d} k^{(1)}_j$.  Since $K_1 \ll \langle k^{(1)}\rangle^{1/(1-\alpha)}$, we have
$$
|k^{(1)}_1| \sim |k^{(2)}_1|, \ \ or \ \  |k^{(1)}_1| \sim |k^{(3)}_1|.
$$

{\it Case 1.} We consider the case $|k^{(1)}_1| = \max_{1\leq i\leq d} |k^{(1)}_i| \sim |k^{(2)}_1|$.

{\it Case 1.1.}  $ |k^{(1)}_1| = \max_{1\leq i\leq d} |k^{(1)}_i| \sim |k^{(2)}_1| \gg |k^{(3)}_1|$.  One can proceeds in the same way as in
Case 1.1 in Section \ref{EstimateL2} to obtain the result, the details are omitted.

{\it Case 1.2.}  $ |k^{(1)}_1| = \max_{1\leq i\leq d} |k^{(1)}_i| \sim |k^{(2)}_1| \sim |k^{(3)}_1|$. First, we conclude that $k^{(1)}_1, \ k^{(2)}_1, \ k^{(3)}_1$ cannot have the same signs, i.e., we must have one of the alternative cases as in \eqref{+--}, \eqref{-+-} and \eqref{--+}.

If  \eqref{+--} is satisfied, by H\"older's inequality, and \eqref{+--2}, bilinear estimates, we have

\begin{align}
& \int_{\mathbb{R}^{d+1}}  \prod^{3}_{\ell=1} \Box^\alpha_{k^{(\ell)}} u_{\ell} \overline{\Box^{\widetilde{\alpha}}_{k} v} dxdt \nonumber\\
& \leq \|\overline{\Box^{\widetilde{\alpha}}_{k} v} \Box^\alpha_{k^{(2)}} \overline{u}  \|_{L^2_{x,t}} \|\Box^\alpha_{k^{(1)}} u   \Box^\alpha_{k^{(3)}} u  \|_{L^2_{x,t}}    \nonumber\\
& \lesssim \ln \langle k^{(1)}\rangle \ \langle k^{(1)} \rangle^{-1/(1-\alpha)}  \langle k  \rangle^{(d-1)\widetilde{\alpha}/2(1-\widetilde{\alpha})} \langle k^{(3)} \rangle^{(d-1)\alpha/2(1-\alpha)} \| \Box^{\widetilde{\alpha}}_{k} v  \|_{V^2_\Delta} \|\Box^\alpha_{k^{(3)}} u  \|_{U^2_\Delta} \|\Box^\alpha_{k^{(1)}} u  \|_{U^2_\Delta} \| \Box^\alpha_{-k^{(2)}} u  \|_{U^2_\Delta}  \nonumber\\
& \lesssim \ln \langle k^{(1)}\rangle \ \langle k^{(1)} \rangle^{(-1-2s)/(1-\alpha)}  \langle k  \rangle^{((d-1)\widetilde{\alpha}/2 +s+ d(\widetilde{\alpha}-\alpha)/2)/(1-\widetilde{\alpha})}   \langle  k  \rangle^{(-s-d(\widetilde{\alpha}-\alpha)/2) /(1-\widetilde{\alpha})}\| \Box^\alpha_{k} v  \|_{V^2_\Delta}   \nonumber\\
& \quad \times   \langle k^{(3)} \rangle^{(  d\alpha/2-s -\alpha/2)/(1-\alpha)} \prod^{3}_{\ell=1} \langle k^{(\ell)} \rangle^{s/(1-\alpha)} \|\Box^\alpha_{\pm k^{(\ell)}} u_{\ell}\|_{U^2_\Delta }. \label{nonlinearkappa=1L2}
\end{align}
Using the same way as in \eqref{L2est}, we have
\begin{align}
\mathscr{L}_2 (u,v)
& \lesssim   \|   v  \|_{Y^{-s-d(\widetilde{\alpha}-\alpha)/2,\widetilde{\alpha}}_\Delta}   \| u \|^{3}_{X^{s,\alpha}_\Delta}. \label{L2estkappa=1}
\end{align}

\subsection{Proof of Theorem \ref{Theorem1}}

Now the proof of Theorem \ref{Theorem1} seems quite standard after the multi-linear estimate is established. We consider the mapping
$$
\mathcal{T}: u(t) \to  S(t)u_0 +{\rm i}\mathscr{A} (|u|^{2\kappa}u)
$$
in the space
$$
\mathcal{D}=\{ u\in X^{s,\alpha}_\Delta ([0,\infty)) : \|u\|_{ X^{s,\alpha}_\Delta ([0,\infty)) } \leq M\}.
$$
Using the same way as in the proof of Theorem 1.1 in \cite{HaHeKo09}, we can show that $\mathcal{T}$ is a contraction mapping on $\mathcal{T}$ to obtain the global existence of solution in $X^{s,\alpha}_\Delta ([0,\infty))$.
To prove the uniqueness, we need the following (see \cite{HaHeKo09})

\begin{lem}[\rm  Right continuity]. \label{rightc}
 Let $T>0$, $u(0)=0$ and $u\in X^{s,\alpha}_{\Delta}([0,T])$. For any $\varepsilon>0$, there exists $0<t_0<T$ such that $\|u\|_{X^{s,\alpha}_{\Delta}([0,t_0])} \leq \varepsilon $.
\end{lem}
Using Lemma \ref{rightc} and following \cite{HaHeKo09}, we see that for two solutions $u,v \in X^{s,\alpha}_\Delta ([0,\infty))$ of NLS, we have
$$
u-v = {\rm i}\mathscr{A} (|u|^{2\kappa}u - |v|^{2\kappa}v).
$$
Applying Lemmas \ref{rightc} and \ref{Nonlinear},  one has that
$$
\|u-v\|_{X^{s,\alpha}_{\Delta}([0,t_0])}   \lesssim  (\|u\|^{2\kappa}_{X^{s,\alpha}_{\Delta}([0,t_0])} + \|v\|^{2\kappa}_{X^{s,\alpha}_{\Delta}([0,t_0])})   \|u-v\|_{X^{s,\alpha}_{\Delta}([0,t_0])} \lesssim  \varepsilon^{2\kappa } \|u-v\|_{X^{s,\alpha}_{\Delta}([0,t_0])}.
$$
This implies the uniqueness of $u\in X^{s,\alpha}_{\Delta}$.

\section{Proof of Theorem \ref{Corollary1}} \label{pfcor}

If $\kappa >2/d$ and $0\leq \alpha<1$, then we have $s_\kappa < s(\kappa)$.  Let $s = s_\kappa+$.  Now we construct a suitable initial value. Choose a sequence $\{k_j\}$ satisfying $k_j \in \mathbb{N}$ and $\langle k_j\rangle^{\alpha/(1-\alpha)} k_j \in [2^{j+1/4}, 2^{j+1/2})$.  Denote
$$
\mathbb{K}^d = \{(k_j,0,...,0): \ j\in J \mathbb{N}\},
$$
where $J$ will be chosen as below.
Put for some small constant $0<\varepsilon\ll 1$,
$$
\widehat{u}_0 = \varepsilon   \sum_{k\in \mathbb{K}^d} \frac{1}{\ln^{2} |k|} \langle k\rangle^{-(s+d\alpha/2)/(1-\alpha)} \mathscr{F}^{-1} \rho^\alpha_k, $$
where $\rho^\alpha_k =\varphi((\cdot- \langle k\rangle^{ \alpha /(1-\alpha)}k)/ c\langle k\rangle^{ \alpha /(1-\alpha)})$, $\varphi$ is as in \eqref{phialpha}  and $c=1/8$.  We can write $u_0$  as
$$
u_0(x) = \varepsilon c^d  \sum_{k\in \mathbb{K}^d} \frac{1}{\ln^{2} |k|} \langle k\rangle^{(-s+d\alpha/2)/(1-\alpha)} e^{i x \langle k\rangle^{ \alpha /(1-\alpha)}k }(\mathscr{F}^{-1} \varphi)(c \langle k\rangle^{ \alpha /(1-\alpha)}  x).
$$
We may assume that $(\mathscr{F}^{-1} \varphi)(0)=1$. We easily see that  for $\alpha>0$
$$
u_0(0) =  \varepsilon   c^d  \sum_{k\in \mathbb{K}^d} \frac{1}{\ln^{2} |k|} \langle k\rangle^{(-s+d\alpha/2)/(1-\alpha)} =\infty
$$
and for $\alpha=0$
$$
|u_0(x)|  \leq   \varepsilon   c^d  \sum_{k\in \mathbb{K}^d} \frac{1}{\ln^{2} |k|}   <\infty.
$$
We see that for $\Lambda=\{\ell\in \mathbb{Z}^d: \ |\ell|\leq C\}$,
$$
\|u_0\|_{M^s_{2,1}} \lesssim  \varepsilon   \sum_{k\in \mathbb{K}^d+ \Lambda} \frac{1}{\ln^{2} |k|}  \langle k\rangle^{- d\alpha/2(1-\alpha)} \|\rho^\alpha_k\|_2 \sim \varepsilon \ll 1.
$$
Moreover, we see that for any $s'>s$ and for $k\in \mathbb{K}^d$,  $|k| \gg 1$,
$$
\|u_0\|_{B^{s'}_{2,\infty}}  \gtrsim  \varepsilon     \frac{1}{\ln^{2} |k|}  \langle k\rangle^{(s'-s)/ (1-\alpha)}   \langle k\rangle^{- d\alpha/2(1-\alpha)} \|\rho^\alpha_k\|_2  \gtrsim   \frac{\varepsilon }{\ln^{2} |k|}  \langle k\rangle^{(s'-s)/ (1-\alpha)}.
$$
It follows that
$$
\|u_0\|_{B^{s'}_{2,\infty}}  =\infty.
$$
In view of Theorem \ref{Theorem1}, one sees that \eqref{NLSi} has a global solution $u\in C(\mathbb{R}, M^{s,\alpha}_{2,1})$.
  Considering the scaling solution $u_\sigma = \sigma^{1/\kappa} u(\sigma^2 t, \sigma x)$, which solves nonlinear
Schr\"odinger equation (NLS):
\begin{align}
 {\rm i}v_t +\Delta v+ \lambda |v|^{2\kappa} v =0,  \quad v(0)= \sigma^{1/\kappa} u_0(\sigma \ \cdot).  \label{NLSs}
\end{align}
We see that
$$
v(0,x) = \varepsilon  c^d  \sigma^{1/\kappa} \sum_{k\in \mathbb{K}^d} \frac{1}{\ln^{2} |k|} \langle k\rangle^{(-s+d\alpha/2)/(1-\alpha)} e^{i x \sigma \langle k\rangle^{ \alpha /(1-\alpha)}k }(\mathscr{F}^{-1} \varphi)( c \sigma \langle k\rangle^{ \alpha /(1-\alpha)}  x).
$$
One can rewrite $v(0, \cdot)$ as
$$
\widehat{v}(0,\xi) = \varepsilon  \sigma^{1/\kappa -d} \sum_{k\in \mathbb{K}^d} \frac{1}{\ln^{2} |k|} \langle k\rangle^{(-s-d\alpha/2)/(1-\alpha)} \varphi \left(\frac{\xi- \sigma \langle k\rangle^{ \alpha /(1-\alpha)}k }{c \sigma \langle k\rangle^{ \alpha /(1-\alpha)}}   \right).
$$

Let $M\gg 1$  and assume without loss of generality that $\sigma =M/\varepsilon^{\kappa/(1-\alpha)} = 2^J$. We easily see that for $\alpha=0$,
\begin{align}
|v (0,0)| \geq  c^d (M^{1/\kappa}/\varepsilon ) |u_0(0)|  \geq c^d M^{1-/\kappa}, \label{large}
\end{align}
and for $\alpha>0$,
\begin{align}
 v (0,0)  =\infty.  \label{large2}
\end{align}
Denote
\begin{align}
\Lambda_k = \left\{l: \ \eta^{\alpha}_l (\xi)  \varphi \left(\frac{\xi- \sigma \langle k\rangle^{ \alpha /(1-\alpha)}k }{ c \sigma \langle k\rangle^{ \alpha /(1-\alpha)}}   \right)   \not\equiv 0 \right\} .  \label{largem}
\end{align}
It is easy to see that for $k\in \mathbb{K}^d$,
\begin{align}
\# \Lambda_k \sim  \sigma^{d(1-\alpha)};  \ \ \langle l\rangle \sim \sigma^{1-\alpha} \langle k\rangle \ \ for \  \ l\in \Lambda_k.  \label{largem2}
\end{align}
We have
\begin{align}
\|v(0,\cdot)\|_{M^{s,\alpha}_{2,1}} & \geq \varepsilon  \sigma^{1/\kappa -d} \ln^{-2}|k| \langle k\rangle^{-(s+d\alpha/2)/(1-\alpha)} \sum_{l\in \Lambda_k} \langle l\rangle^{s/(1-\alpha)} \left\|\eta^{\alpha}_l  \varphi \left(\frac{\xi- \sigma \langle k\rangle^{ \alpha /(1-\alpha)}k }{c \sigma \langle k\rangle^{ \alpha /(1-\alpha)}}   \right) \right\|_2 \nonumber\\
& \gtrsim \varepsilon  \sigma^{1/\kappa+s - d\alpha/2} \ln^{-2}|k|    \geq \varepsilon  \sigma^{(1-\alpha)+/\kappa}  \ln^{-2}|k|.
\end{align}
Taking $|k|^{1/(1-\alpha)} \sim 2$, we immediately have
\begin{align}
\|v(0,\cdot)\|_{M^{s,\alpha}_{2,1}}  \geq  M^{(1-\alpha)/\kappa}.
\end{align}

\section{Ill Posedness}

Now let $u(\delta,t)$ satisfy
\begin{align}
& {\rm i}u_t +\Delta u+ \lambda |u|^{2\kappa} u =0,  \quad u(0)=\delta v_0. \label{deltaNLS}
\end{align}
One has that
\begin{align}
\frac{\partial^{2\kappa+1}u(0,t)}{\partial \delta^{2\kappa+1}}  = \int^t_0 S(t-\tau) |S(\tau) v_0|^{2\kappa} S(\tau) v_0 d\tau.  \label{derivative}
\end{align}
It follows that
\begin{align}
\widehat{\frac{\partial^{2\kappa+1}u(0,t)}{\partial \delta^{2\kappa+1}}}  = e^{-{\rm i}t|\xi|^2} \int_{\mathbb{R}^{2\kappa d}} \frac{e^{-{\rm i}tP}-1}{P} \widehat{v}_0(\xi-\xi_1-...-\xi_{2\kappa}) \prod^{2\kappa}_{j=1} \widehat{v}_0(\xi_j) d\xi_1...d\xi_{2\kappa},  \label{derivativeF}
\end{align}
where
$$
P= -\sum^\kappa_{j=1} |\xi_j|^2 -|\xi-\xi_1-...-\xi_{2\kappa}|^2 + \sum^{2\kappa}_{j=\kappa+1} |\xi_j|^2 + |\xi|^2.
$$
Take $N\gg 1$, $k=(N,N,...,N)$  and
\begin{align}
\widehat{v}_0 (\xi)  = \langle k\rangle^{-(s+d\alpha/2)/(1-\alpha)} (\varphi_k (\xi) + \varphi_{-k}(\xi)),
\end{align}
where we denote by $\varphi$ a smooth cut-off function adapted to the unit ball as in \eqref{varphi} and
$$
\varphi_k (\xi) =\varphi \left(\frac{\xi- \langle k\rangle^{\alpha/(1-\alpha)} k }{ c \langle k\rangle^{\alpha/(1-\alpha)}}\right).
$$

By the definition of norm on $M^{s,\alpha}_{2,1}$, it is easy to see that
$$
\|\delta v_0\|_{M^{s,\alpha}_{2,1}} \sim \delta.
$$
Again, in view of the definition of norm on $M^{s,\alpha}_{2,1}$, we have
\begin{align}
&  \left\| \frac{\partial^{2\kappa+1}u(0,t)}{\partial \delta^{2\kappa+1}}  \right\|_{M^{s,\alpha}_{2,1}}
\gtrsim   \langle k\rangle^{(s-(2\kappa+1)(s+d\alpha/2)) /(1-\alpha)} \nonumber\\
&  \times \left\| \varphi_k (\xi)     \int_{\mathbb{R}^{2\kappa d}} \frac{e^{-{\rm i}tP}-1}{P} \varphi_k(\xi-\xi_1-...-\xi_{2\kappa}) \prod^{ \kappa}_{j=1} \varphi_{-k}(\xi_j)   \prod^{ 2\kappa}_{j=\kappa+1} \varphi_k(\xi_j) d\xi_1...d\xi_{2\kappa}  \right\|_{2}. \label{derivativeFN}
\end{align}
By change of variables
\begin{align}
& \xi_j +  \langle k\rangle^{\alpha/(1-\alpha)}k = c \langle k\rangle^{\alpha/(1-\alpha)} \zeta_j , \ j=1,...,\kappa,  \nonumber\\
& \xi_j -  \langle k\rangle^{\alpha/(1-\alpha)}k = c \langle k\rangle^{\alpha/(1-\alpha)} \zeta_j , \ j=\kappa+1,...,2\kappa,  \nonumber \\
& \xi  -  \langle k\rangle^{\alpha/(1-\alpha)}k = c \langle k\rangle^{\alpha/(1-\alpha)} \zeta,    \nonumber
\end{align}
we see that
$$
P= \langle k\rangle^{2\alpha/(1-\alpha)} \left(- \sum^{\kappa}_{j=1} |\zeta_j|^2 + \sum^{2\kappa}_{j=\kappa+1} |\zeta_j|^2 - |\zeta-\zeta_1-...-\zeta_{2\kappa}|^2 + |\zeta|^2 \right).
$$
By choosing $t\sim \langle k\rangle^{- 2\alpha/(1-\alpha)}$, from \eqref{derivativeFN} it follows that
\begin{align}
\left\| \frac{\partial^{2\kappa+1}u(0,t)}{\partial \delta^{2\kappa+1}}  \right\|_{M^{s,\alpha}_{2,1}}
& \gtrsim  \langle k\rangle^{(s-(2\kappa+1)(s+d\alpha/2) /(1-\alpha)} \langle k\rangle^{ 2\kappa d\alpha /(1-\alpha)} \langle k\rangle^{ d\alpha /2(1-\alpha)}  \langle k\rangle^{- 2\alpha/(1-\alpha)} \nonumber\\
& \ \  \times \left\| \eta  (\zeta)     \int_{\mathbb{R}^{2\kappa d}}   \eta (\zeta-\zeta_1-...-\zeta_{2\kappa}) \prod^{ 2\kappa}_{j=1} \eta (\eta_j)     d\zeta_1...d\zeta_{2\kappa}  \right\|_{2} \nonumber\\
& \gtrsim  \langle k\rangle^{ [2\kappa (d\alpha/2-s) -2\alpha] /(1-\alpha)}.    \label{derivativeFNE}
\end{align}
Hence, if $s<d\alpha/2-\alpha/\kappa$, we have for $N\gg 1$,
\begin{align}
\left\| \frac{\partial^{2\kappa+1}u(0,t)}{\partial \delta^{2\kappa+1}}  \right\|_{M^{s,\alpha}_{2,1}}
\geq N^{\sigma}, \ \ {\rm for \ some} \ \  \sigma>0.
\end{align}
In particular, in the case $\alpha=0$, we have
\begin{align}
\left\| \frac{\partial^{2\kappa+1}u(0,t)}{\partial \delta^{2\kappa+1}}  \right\|_{M^{s,\alpha}_{2,1}}
 \gtrsim  \langle k\rangle^{-2\kappa s}.    \label{derivativeFNE0}
\end{align}
Noticing that $\|v_0\|_{M^{-s/2}_{2,1}} \sim N^{s/2} \ll 1$ if $s < s_{\kappa}=0$ (we can assume that $s=0+$),  we have NLS is globally well-posed in $C(\mathbb{R}, M^{-s/2}_{2,1})$.  So we have a Taylor expansion of $u(\delta,t)$ at $\delta=0$:
$$
u(\delta,t)= \delta S(t)v_0 + \sum^\infty_{j=2\kappa+1}  \frac{\partial^{j}u (0,t)}{\partial \delta^{j}}  \delta^j, \ \ {\rm in} \ \  C(\mathbb{R}, M^{-s/2}_{2,1}).
$$
Using Theorem \ref{Theorem1}, we have
$$
\left\|u(\delta,t) - \frac{\partial^{2\kappa+1}u(0,t)}{\partial \delta^{2\kappa+1}} \right\|_{M^{-s/2}_{2,1}}  \leq \|u\|_{M^{-s/2}_{2,1}}  +  \left\|\frac{\partial^{2\kappa+1}u(0,t)}{\partial \delta^{2\kappa+1}} \right\|_{ M^{-s/2}_{2,1}} \delta^{2\kappa+1}.
$$
Applying the multi-linear estimate,
$$
\left\|\frac{\partial^{2\kappa+1}u(0,t)}{\partial \delta^{2\kappa+1}} \right\|_{ M^{-s/2}_{2,1}}  \leq C \|S(t)v_0\|^{2\kappa+1}_{X^{-s/2,0}_{\Delta}} \lesssim 1.
$$
However, we have for some fixed $t\ll 1$,
\begin{align}
\|u(\delta, t)\|_{M^s_{2,1}} & \geq   \left\|\frac{\partial^{2\kappa+1}u(0,t)}{\partial \delta^{2\kappa+1}} \right\|_{ M^{s}_{2,1}} \delta^{2\kappa+1} - \left\|u(\delta,t) - \frac{\partial^{2\kappa+1}u(0,t)}{\partial \delta^{2\kappa+1}} \delta^{2\kappa+1} \right\|_{M^{-s/2}_{2,1}}  \nonumber\\
& \geq  N^{-2\kappa s}   \delta^{2\kappa+1} - C.
\end{align}
Since $\delta$ is independent of $N$, we have  $
\|u(\delta, t)\|_{M^s_{2,1}} \gg 1$ by choosing  $ N^{-2\kappa s} \gg  \delta^{-2\kappa+1}$. This implies that the  map $\delta \to u(\delta, t)$ is discontinuous at $\delta=0$. \\

{\bf Acknowledgement. \rm  B. Wang is supported in part by an NSFC grant 11771024. Part work of this paper has been talked in the ``Workshop on nonlinear evolution equations" which was held in Kumamoto University during January 20-23, 2018, B. Wang is very grateful to Professor N. Hayashi for his invitation and support in the conference. }

\end{document}